\documentclass[amssymb,twoside,12pt,reqno,amscd]{amsart}
\usepackage{amsmath}
\usepackage{amscd}
\usepackage{graphics}
\usepackage{latexsym}

\theoremstyle{plain}
\newtheorem{theorem}{Theorem}[section]
\newtheorem{thm}[theorem]{Theorem}
\newtheorem{cor}[theorem]{Corollary}
\newtheorem{lem}[theorem]{Lemma}
\newtheorem{prop}[theorem]{Proposition}
\newtheorem{defn}[theorem]{Definition}

\theoremstyle{definition}

\theoremstyle{remark}

\newcommand{\ZZ}{\mathbb{Z}}
\newcommand{\CC}{\mathbb{C}}

\newcommand{\AAA}{\mathbb{A}}
\newcommand{\PP}{\mathbb{P}}

\newcommand{\SP}{\text{Spec }}

\newcommand{\Mgnb}[2]{\overline{\mathcal{M}}_{#1,#2}}

\newcommand{\Kgnb}[4]{\overline{\mathcal{M}}_{#1,#2}(#3,#4)}

\newcommand{\Kbm}[2]{\overline{\mathcal{M}}(#1,#2)}

\newcommand{\Kbmo}[2]{{\mathcal{M}}(#1,#2)}

\newcommand{\lt}{\left}
\newcommand{\rt}{\right}

\newcommand{\mc}{\mathcal}
\newcommand{\cB}{\mathcal B}
\newcommand{\cE}{{\mathcal E}}

\newcommand{\cX}{{\mathcal X}}
\newcommand{\cD}{{\mathcal D}}
\newcommand{\lci}{\mathcal {LCI}}
\newcommand{\fe}{\mathcal {FE}}

\newcommand{\OO}{\mathcal O}

\begin{document}

\title{Rational Curves On Hypersurfaces of Low Degree}
\author[Harris]{Joe Harris} \address{Department of Mathematics \\ 
  Harvard University \\ Cambridge MA 02138}
\email{harris@math.harvard.edu} \author[Roth]{Mike Roth}
\address{Department of Mathematics \\ University of Michigan \\ Ann
  Arbor, MI 48109} \email{mikeroth@math.lsa.umich.edu}
\author[Starr]{Jason Starr}\thanks{The third author was partially
  supported by an NSF Graduate Research Fellowship and a Sloan
  Dissertation Fellowship} \address{Department of Mathematics \\ 
  Massachusetts Institute of Technology \\ Cambridge MA 02139}
\email{jstarr@math.mit.edu} \date{\today}

\begin{abstract}
  We prove that for $n>2$ and $d<\frac{n+1}{2}$, a \emph{general}
  complex hypersurface $X\subset \PP^n$ of degree $d$ has the property
  that for each integer $e$ the scheme $R_e(X)$ parametrizing degree
  $e$, smooth rational curves on $X$ is an integral, local complete
  intersection scheme of ``expected'' dimension $(n+1-d)e+(n-4)$.

\

The techniques used in the proof include:
\begin{enumerate}
\item  Classical results about lines on hypersurfaces including a new
result about flatness of the projection map from the space of
\emph{pointed lines}.
\item  The Kontsevich moduli space of stable maps,
$\Kgnb{0}{r}{X}{e}$.  In particular we use the
deformation theory of stable maps, properness of the stack
$\Kgnb{0}{r}{X}{e}$, and the decomposition of
$\Kgnb{0}{r}{X}{e}$ described in ~\cite{BM}.
\item  A version of Mori's bend-and-break lemma.
\end{enumerate}
\end{abstract}

\maketitle

\section{Summary}~\label{sec-sum}

\subsection{Brief Summary}~\label{subsec-brsum}
All schemes we consider will be $\CC$-schemes and all morphisms will
be morphisms of $\CC$-schemes.  All (absolute) products will be over
$\CC$.

\

For a projective scheme $X$ over $\CC$ along with an ample line bundle
$L$ we define $R_e(X)$ to be the open subscheme of the Hilbert scheme
$\text{Hilb}^{et+1}(X/k)$ which parametrizes smooth rational curves of
degree $e$ lying in $X$.

\begin{thm}~\label{thm-main}  Let $n > 2$ be an integer and let $d$ be a
  positive integer such that $d < \frac{n+1}{2}$.  For a general
  hypersurface $X\subset \PP^n$ of degree d and for every integer
  $e\geq 1$, the scheme $R_e(X)$ is an integral, local complete
  intersection scheme of dimension $(n+1-d)e+(n-4)$.
\end{thm}

The idea of the proof is as follows.  There is an embedding of
$R_e(X)$ into the smooth scheme $R_e(\PP^n)$. Denote by
$\pi:U_e(\PP^n)\rightarrow R_e(\PP^n)$ the universal family of
rational curves in $\PP^n$ and by $\rho:U_e(\PP^n)\rightarrow \PP^n$
the \emph{evaluation} morphism.  Then $R_e(X)$ is the scheme of zeroes
of a section of the locally free sheaf $\pi_* \rho^* \OO_{\PP^n}(d)$.
Thus to prove that $R_e(X)$ is a local complete intersection scheme,
it suffices to prove that the codimension of $R_e(X)$ in $R_e(\PP^n)$
equals the rank of $\pi_* \rho^* \OO_{\PP^n}(d)$.

\

The remainder of the proof is a ``deformation and specialization''
argument: we embed the non-proper scheme $R_e(X)$ as an open subscheme
of a proper scheme which is still modular, i.e. we choose a ``modular
compactification''. Then we show that every generic point in $R_e(X)$
specializes to a point in the ``boundary'' of the compactification.
We use deformation theory to study the irreducible components of the
boundary of the compactification.  In particular we show that a
general point of each irreducible component of the boundary is a
unibranch point of the compactification whose local ring is reduced
and has the expected dimension.  This reduces the proof to a
combinatorial argument.

\subsection{Detailed Summary}~\label{subsec-detsum}

In the next few paragraphs we will give a detailed summary of the
proof.  Our compactification consists of the embedding of $R_e(X)$ as
an open subscheme in the Kontsevich moduli space $\Kgnb{0}{0}{X}{e}$
parametrizing stable maps to $X$.  We recall the partition of
$\Kgnb{0}{0}{X}{e}$ into locally closed subsets defined in~\cite{BM};
we call this partition the \emph{Behrend-Manin decomposition} (our
partition differs slightly from that in~\cite{BM}).  In particular,
the image of $R_e(X)$ is a dense open subset of a component of this
partition. We identify certain \emph{basic components} as those
components of the partition parametrizing stable maps such that each
irreducible component of the domain curve is mapped to a line in $X$.

\

We prove a new result about lines on $X$.  We define the incidence
correspondence of \emph{pointed lines} in X:
\begin{equation}
F_{0,1}(X)=\{(p,l) |p \text{ a point}, l \text{ a line}, p\in
l\subset X\}.
\end{equation}
We prove that for a general hypersurface $X\subset \PP^n$ of degree
$d\leq n-1$, the projection morphism $F_{0,1}(X)\rightarrow X$ is flat
of relative dimension $n-d-1$.  From this theorem it easily follows
that each basic component $B$ is an integral scheme whose general
point is a unibranch point of $\Kgnb{0}{0}{X}{e}$ at which
$\Kgnb{0}{0}{X}{e}$ is reduced of dimension $(n+1-d)e+(n-4)$. Thus for
each basic component $B$ there is a unique irreducible component
$M(B)$ of $\Kgnb{0}{0}{X}{e}$ which contains $B$, and $M(B)$ is
reduced and has dimension $(n+1-d)e+(n-4)$.

\

Using a version of the bend-and-break lemma of Mori, we prove that
every irreducible component of $\Kgnb{0}{0}{X}{e}$ is of the form
$M(B)$ for some basic component $B$.  Using this fact and results
about flatness, we bootstrap to prove that each evaluation map
$\Kgnb{0}{r}{X}{e}\rightarrow X$ is flat of the expected dimension and
is \emph{generically unobstructed}.  This implies that each
$\Kgnb{0}{r}{X}{e}$ (including $r=0$) has the expected dimension and
is generically smooth.  Thus $\Kgnb{0}{r}{X}{e}$ is a reduced, local
complete intersection stack of the expected dimension and it only
remains to prove that $\Kgnb{0}{r}{X}{e}$ is irreducible, i.e. it
remains to prove that all of the irreducible components $M(B)$ are
actually equal.

\

To prove that all of the irreducible components $M(B)$ are equal, we
observe that there is a combinatorially defined equivalence relation
defined on the set of basic components $B$ such that equivalent basic
components, $B\cong B'$ satisfy $M(B)=M(B')$.  Thus we are reduced to
a combinatorial argument which proves that all basic components are
equivalent.

\

Along the way we generalize the strategy of proof above so that it
could apply to smooth projective schemes $X$ other than hypersurfaces
$X\subset \PP^n$ of degree $d< \frac{n+1}{2}$ (this is made completely
explicit for complete intersections in $\PP^n$).  One is reduced to
proving
\begin{enumerate}
\item The evaluation morphism $\Kgnb{0}{1}{X}{e}\rightarrow X$ is
  flat, generically unobstructed and the general fiber is
  geometrically irreducible.
\item For each positive integer $e$ at most the \emph{threshold
    degree} of $X$, the evaluation morphism
    $\Kgnb{0}{1}{X}{e}\rightarrow X$ is flat of the expected
    dimension.
\item For each positive integer $e$ at most the threshold degree of
    $X$, the stack $\Kgnb{0}{0}{X}{e}$ is irreducible.
\end{enumerate}
The most difficult condition to verify seems to be $(2)$, but it is
our hope that this can be verified for a larger class of Fano schemes
$X$ than the hypersurfaces above.

\

This paper can be seen as a partial generalization of the
irreducibilty results ~\cite[corollary 1]{KP} to hypersurfaces
$X\subset \PP^n$ of degree roughly $d\leq\frac{n}{2}$. In a
forthcoming paper, \cite{HS2}, we give a partial generalization of the
rationality results of~\cite[theorem 3]{KP} to hypersurfaces $X\subset
\PP^n$ of degree roughly $d\leq \sqrt{n}$.

\subsection{Notation}~\label{subsec-not}
Given a $\CC$-vector space $W$, $\PP W$ denotes the projective space
$\text{Proj}(\oplus_{d\geq 0} S^d(W^*))$ which parametrizes
one-dimensional linear subspaces of $W$ (not one-dimensional quotient
spaces of $W$).  Given any integers $k\leq n$, $G(k,n)$ denotes the
Grassmannian which parametrizes $k$-dimensional linear subspaces of
$\CC^n$.  For a triple of integers $k\leq l\leq n$, $F((k,l),n)$
denotes the partial flag variety which parametrizes partial flags
$V_1\subset V_2\subset \CC^n$ of linear subspaces with
$\text{dim}(V_1)=k$ and $\text{dim}(V_2)=l$.

\subsection{Acknowledgments}~\label{subsec-acknow}  
The authors would like to thank Tom Graber, Ravi Vakil, Andreas
Gathmann, Olivier Debarre, and especially Johan de Jong for many
useful discussions.

\section{Lines on Hypersurfaces}~\label{sec-lines}

We denote $W=H^0(\PP^n,\OO_{\PP^n}(d))$ so that $\PP W$ parametrizes
degree $d$ hypersurfaces $X\subset \PP^n$.  Let $\cX \subset \PP W
\times \PP^n$ denote the universal family of degree $d$ hypersurfaces
in $\PP^n$.  For each degree $d$ hypersurface $X\subset \PP^n$, denote
by $F_1(X)$ the subscheme of $G(2,n+1)$ which parametrizes lines
$L\subset X\subset \PP^n$.  Denote by $F((1,2),n+1)$ the partial flag
variety parametrizing pairs $(p,L)$ where $L\subset \PP^n$ is a line
and $p\in L$ is a point.  Denote by $F_{0,1}(X)$ the subscheme of
$F((1,2),n+1)$ which parametrizes pairs $(p,L)$ with $p\in L \subset
X\subset \PP^n$.  

\

Similarly, let $F_1(\cX)\subset \PP W\times
G(2,n+1)$ denote the subscheme parametrizing pairs $(X,L)$ with $L\in
F_1(X)$, and let $F_{0,1}(\cX)\subset \PP W\times F((1,2),n+1)$ denote
the subscheme parametrizing triples $(X,p,L)$ with $(p,L)\in
F_{0,1}(X)$.  There are projection morphisms
\begin{eqnarray}
\pi_0:\PP W \times F((1,2),n+1)\rightarrow \PP W, \\
\pi_1:\PP W \times F((1,2),n+1)\rightarrow \PP^n, \\
\pi_2:\PP W \times F((1,2),n+1)\rightarrow G(2,n+1)
\end{eqnarray}
By construction, the morphism $(\pi_0,\pi_1):F_{0,1}(\cX)\rightarrow
\PP W\times \PP^n$ factors through $\cX \subset \PP W\times \PP^n$.
Denote by $\rho:F_{0,1}(\cX)\rightarrow \cX$ the induced morphism.
For a particular hypersurface $X\in \PP W$, denote by
$\rho_X:F_{0,1}(X)\rightarrow X$ the fiber of $\rho$.

\

The main result of this section is the following theorem:
\begin{thm} \label{thm-thm1}  Let $d$ be a positive integer with $d \leq
n-1$.  For a 
general hypersurface $X\in \PP W$ the morphism
$\rho_X:F_{0,1}(X)\rightarrow X$ is flat of relative dimension
$n-d-1$.
\end{thm}

We give the proof in the remainder of this section.  From now on
we assume that $d$ is given such that $d\leq n-1$.

\

Denote by $\OO$ the structure sheaf of $\PP W \times F((1,2),n+1)$ and
denote by $\OO_F$ the $\OO$-module which is the pushforward of the
structure sheaf of $F_{0,1}(\cX)$.  On $G(2,n+1)$ we have a universal
rank 2 subbundle of $\OO_{G(2,n+1)}^{\oplus n+1}$.  Denote by $S$ the
pullback under $\pi_2$ of this universal subbundle to $\PP W\times
F((1,2),n+1)$.  And denote by $U$ the pullback under $\pi_0$ of the
universal rank 1 subbundle $\OO_{\PP W}(-1) \subset
H^0(\PP^n,\OO_{\PP^n}(d))\otimes_\CC \OO_{\PP W}$.  By restricting a
section of $H^0(\PP^n,\OO_{\PP^n}(d))$ to a line $L\subset \PP^n$
parametrized by a point $[L]\in G(2,n+1)$, we get a map $U\rightarrow
\text{Sym}^d(S^\vee)$.  By adjunction, this gives rise to a map
$\text{Sym}^d(S)\otimes_\OO U \rightarrow \OO$, whose image is exactly the
ideal sheaf of $F_{0,1}(\cX)$.  In other words, there is a partial
resolution of coherent sheaves:
\begin{equation}
\begin{CD}
  \text{Sym}^d(S)\otimes_{\OO} U @> \sigma >> \OO @>>> \OO_F @>>> 0.
\end{CD}
\end{equation}
In other words, $F_{0,1}(\cX)$ is the zero scheme of the global
section $\OO\rightarrow \text{Sym}^d(S^\vee)\otimes_\OO U^\vee$ which is
the transpose of $\sigma$.

\

Since $\text{Sym}^d(S)\otimes_\OO U$ is locally free of rank $d+1$, every
irreducible component of $\cX$ has codimension at most $d+1$ in $\PP
W\times F((1,2),n+1)$.  Therefore, every (nonempty) fiber of $\rho$
has dimension at least $n-d-1$.  We define ${\mc U}\subset \cX$ as a
set to be
\begin{equation}
{\mc U} = \lt\{ (X,p)\in \cX |\ \text{dim}(\rho^{-1}(X,p))\leq n-d-1
\rt\}.
\end{equation}
It follows by upper semicontinuity of the fiber dimension that ${\mc
  U}$ is a Zariski open subset of $\cX$, and we give it the
corresponding structure of open subscheme of $\cX$.  A priori ${\mc
  U}$ might contain points $(X,p)\subset \cX$ for which
$\rho^{-1}(X,p)$ is empty.  But by ~\cite[exercise V.4.6]{K}, it
follows that $\rho$ is surjective (also this exercise rederives the
statement above about dimensions of fibers).

\

Notice that the projection morphism $\pi_0:\cX\rightarrow \PP^n$ is a
projective bundle whose fiber over $p\in \PP^n$ is identified, as a
subscheme of $\PP W$, with the hyperplane parametrizing $X\in \PP W$
with $p\in X$.  In particular, $\cX$ is a smooth $k$-scheme.  Given
the map $\sigma$ above, we can form the Koszul complex of locally free
$\OO$-modules in the usual way.  By ~\cite[theorem 17.4 (iii)(4)]{Ma},
this complex is acyclic over ${\mc U}$.  Therefore the fibers of
$\rho$ over ${\mc U}$, considered as subschemes of the appropriate
fiber of $\pi_1:F((1,2),n+1)\rightarrow \PP^n$, all have equal Hilbert
polynomial.  Since ${\mc U}$ is smooth, it follows from ~\cite[theorem
III.9.9]{H} that $\rho$ is flat over ${\mc U}$.

\

Let ${\mc Y} \subset F_{0,1}(\cX)$ denote the complement of ${\mc U}$
with the induced, reduced scheme structure.  Theorem 2 is equivalent
to the statement that $\pi_0|_{\mc Y}:{\mc Y}\rightarrow \PP W$ is not
surjective.  Denote by $e$ the codimension of $\mathcal{Y}$ in
$\mathcal{X}$.  Since the fiber dimension of $\mathcal{X}\rightarrow
\PP W$ is $n-1$, to prove that $\mathcal{Y}$ fails to dominate $\PP
W$, it suffices to prove that $e > n-1$.  In the remainder of this
section we will prove that $e>n-1$.

\

On $\PP^n$ let $Q$ denote the locally free quotient sheaf of
$\OO_{\PP^n}^{\oplus n+1}$ by $\OO_{\PP V}(-1)$.  The dual injection
$Q^\vee \hookrightarrow \lt(\OO_{\PP^n}^{\oplus n+1}\rt)^\vee$ can be
considered as a filtration of $\lt(\OO_{\PP^n}^{\oplus n+1}\rt)^\vee$.
The $d$th symmetric product of this filtration is a filtration of $W
\otimes_\CC \OO_{\PP^n}$:
\begin{equation}
  W \otimes_\CC \OO = F^{0,d} \supset F^{1,d}
\supset \dots \supset F^{d,d} \supset F^{d+1,d} = 0.
\end{equation}
Here $F^{i,d}$ is the locally free subsheaf of $ W \otimes_\CC \OO$
which is the image of the multiplication map
\begin{equation}
\text{Sym}^{d-i} \lt(\OO_{\PP^n}^{\oplus n+1}\rt)^\vee  \otimes_{\OO_{\PP^n}} 
\text{Sym}^i Q^\vee\rightarrow W\otimes_\CC \OO.
\end{equation}
The associated graded sheaves of this filtration $G^{i,d} =
F^{i,d}/F^{i+1,d}$ are canonically isomorphic to the sheaves
$\OO_{\PP^n}(d-i)\otimes_{\OO_{\PP^n}} \text{Sym}^i Q^\vee$.

\

In particular, notice that $F^{1,d}$ is simply the kernel of the
evaluation map $W\otimes_\CC \OO_{\PP^n}\rightarrow \OO_{\PP^n}(d)$,
i.e. the vector bundle parametrizes pairs $(\Phi,p)$ where $\Phi\in W$
is such that $\Phi(p)=0$.  We identify a nonzero section $\Phi\in W$,
up to nonzero scaling, with the hypersurface it defines $X=V(\Phi)$.
Then the associated projective bundle $\PP F^{1,d}$ inside $\PP
W\times \PP^n$ is the closed subscheme parametrizing pairs $(X,p)$
with $p\in X$, i.e. $\PP F^{1,d} = \cX$.  To prove the inequality
$e>n-1$ from above, it suffices to prove that for each $p\in \PP^n$ 
(equivalently for any $p\in \PP^n$ by homogeneity) the intersection
${\mc Y}\cap \pi_0^{-1}(p) \subset \cX$ has codimension greater than
$n-1$ in $\pi_0^{-1}(p)$.  In the remainder of this section we will
prove this.

\

Observe the filtration above is not split on $\PP^n$.  But we can find
a covering of $\PP^n$ by open affine subschemes $A_\alpha \subset
\PP^n$ over which we do have a splitting (for example, the standard
covering by complements of coordinate hyperplanes).  Here by
\emph{splitting} we mean an isomorphism of bundles over $A_\alpha$
\begin{equation}
\begin{CD} s: W \otimes_\CC \OO_{A_\alpha} @>>>
  \oplus_{j=0}^{d} \OO_{\PP^n}(d-j)\otimes_{\OO_{\PP^n}} \text{Sym}^j
  Q^\vee|_{A_\alpha}
\end{CD}
\end{equation}
which maps $F^{i,d}|_{A_\alpha}$ to the subbundle $\oplus_{j=i}^{d}
\OO_{\PP^n}(d-j)\otimes_{\OO_{\PP^n}} \text{Sym}^j Q^\vee|_{A_\alpha}$
and such that the induced isomorphism $G^{i,d}|_{A_\alpha}\rightarrow
\OO_{\PP^n}(d-i) \otimes_{\OO_{\PP^n}} \text{Sym}^i
Q^\vee|_{A_\alpha}$ is the isomorphism from above.

\

Given an open affine $A_\alpha\subset \PP^n$ we can form the
projective bundle $\PP_{A_\alpha}(F^{1,d}|_{A_\alpha})$ over
$A_\alpha$.  Given a splitting $s$ on $A_\alpha$, denote by
\begin{equation}
\Delta_j(s)\subset \PP_{A_\alpha} (F^{1,d}|_{A_\alpha})
\end{equation}
the closed subscheme which parametrizes pairs $(\Phi,x)$, $x\in
A_\alpha, \Phi\in F^{1,d}|x$ such that the $j$th component of $s(\Phi)$
is zero.  Thus $\Delta_0(s)$ is all of $\PP_{A_\alpha}
(F^{1,d}|_{A_\alpha})$.  And, considering $\PP_{A_\alpha}
(F^{1,d}|_{A_\alpha})$ as an open subscheme of $\mathcal{X}$,
$\Delta_1(s)$ is the intersection of $\PP_{A_\alpha}
(F^{1,d}|_{A_\alpha})$ with the singular locus of the projection
morphism $\pi_0:\mathcal{X}\rightarrow \PP W$.  Although $\Delta_0(s)$
and $\Delta_1(s)$ are independent of the choice of $s$, the same is
not true for $\Delta_i(s)$ with $i>1$.  The next result follows
immediately from the definition of the $\Delta_j(s)$.

\

\begin{lem}\label{lem-lem0}  For $j>0$ the codimension of $\Delta_j(s)$ in
  $\PP_{A_\alpha}(F^{1,d}|_{A_\alpha})$ equals
\begin{equation}
\text{rank } \lt( \OO_{\PP^n}(d-j)\otimes_{\OO_{\PP^n}} \text{Sym}^j Q^\vee
\rt) = \lt( 
\begin{array}{c}  n-1+j \\ n-1 \end{array}
\rt).
\end{equation}
\end{lem}

\

In particular, for $j>0$ the codimension of $\Delta_j(s)$ in
$\PP_{A_\alpha}(F^{1,d}|_{A_\alpha})$ is at least $n-1+j > n-1$.  For
each open subscheme ${A_\alpha}\subset \PP V$, each splitting $s$, and
each point $p\in {A_\alpha}$, define the locally closed subscheme
\begin{equation}
\mathcal{Y}_{p,s}:=\lt( \mathcal{Y}\cap \pi_0^{-1}(p) \rt)-\cup_{j=1}^d
\lt( \Delta_j(s)\cap \pi_0^{-1}(p) \rt).
\end{equation}
To establish that $e>n-1$, it suffices to prove that the codimension
of ${\mc Y}_{p,s}$ as a subscheme of $\pi_0^{-1}(p)$ 
has codimension greater than $n-1$.  In the remainder of this section
we prove this inequality.

\

On the complement of the closed subset $\Delta(s) :=\cup_{j=1}^d \lt(
\Delta_j(s) \rt)$ there is a morphism
\begin{equation}
\PP_{A_\alpha} \lt( F^{1,d}|_{A_\alpha} \rt)-\Delta(s) \xrightarrow{\beta}
  \prod_{j=1}^d \PP_{A_\alpha} \lt( \OO_{\PP^n}(d-j)
  \otimes_{\OO_{\PP^n}} \text{Sym}^j Q^\vee \rt)|_{A_\alpha}.
\end{equation}
We identify the space
\begin{equation}
\PP_{A_\alpha} \lt( \OO_{\PP^n}(d-j) \otimes_{\OO_{\PP^n}} 
   \text{Sym}^j Q^\vee \rt)|_{A_\alpha}
\end{equation}  
with the scheme parametrizing degree $j$ hypersurfaces in fibers of
the projection morphism $\PP_{A_\alpha} Q|{A_\alpha} \rightarrow
{A_\alpha}$.  Thus $\beta$ assigns to each suitable pair $([\Phi],p)$
a sequence of hypersurfaces in $\PP Q|_p$.  Denote this sequence by
$\lt(X_1,\dots,X_j,\dots,X_d\rt)$.  Also, denote by $X=X_1\cap \dots
\cap X_d$ the intersection of these hypersurfaces.

\begin{lem}\label{lem-lem1}  If we denote by $X$ the hypersurface in $\PP V$
  corresponding to $\Phi$, then $X_1\cap \dots\cap X_d$ is the fiber
  of $p\in X$ under $\rho_X:F_{0,1}(X)\rightarrow X$.
\end{lem}

\begin{proof}  
  This is most easily seen by passing to local coordinates.  Let
  $(x_0,\dots,x_n)$ be a system of homogeneous coordinates on $\PP V$
  (i.e. a basis for $V^\vee$) and let $p$ be the point with
  homogeneous coordinates $[0,\dots,0,1]$.  We define a splitting $s$
  as follows: for each degree $d$ homogeneous polynomial $\Phi$ in
  $(x_0,\dots,x_n)$ we have a unique decomposition
\begin{equation}
\Phi = \Phi_d + \Phi_{d-1}x_n+\dots+ \Phi_{d-i}x_n^i +\dots
+\Phi_0 x_n^d
\end{equation}
where each $\Phi_i$ is a homogeneous polynomial of degree $i$ in
$(x_0,\dots,x_{n-1})$.  Then the fiber of $F^{1,d}$ at $p$ consists of
those polynomials such that $\Phi_0=0$ and $\beta(\Phi) =
(\Phi_d,\dots,\Phi_1)$. For any line $L$ passing through $p$ there is
a unique point of the form $y=(a_0,\dots,a_{n-1},0)$ contained in $L$.
Let $\PP^1\rightarrow \PP^n$ be the morphism given by
\begin{equation} 
(t_0,t_1)\mapsto (t_1a_0,t_1a_1,\dots,t_1a_{n-1},t_0+t_1a_n).
\end{equation}  
The image of this morphism is just $L$.  Substituting into $\Phi$
yields the polynomial on $\PP^1$ given by
\begin{eqnarray*}
t_0^d\Phi_d(a_0,\dots,a_{n-1}) + \dots +
t_0^{d-i}t_1^i\Phi_{d-i}(a_0,\dots,a_{n-1})+ \dots \\
+ t_0t_1^{d-1}\Phi_1(a_0,\dots,a_{n-1}).
\end{eqnarray*}   
The line $L$ is contained in $X$ iff this polynomial is identically
zero iff each of the terms $\Phi_i(a_0,\dots,a_n)$ is zero.  One can
show that the homogeneous ideal generated by the terms $\Phi_i$ is
independent of our particular splitting.
\end{proof}

In particular, we conclude that every fiber of $\beta$ which
intersects $\mathcal{Y}$ is contained in $\mathcal{Y}$.  Therefore the
codimension of $\mathcal{Y}_{p,s}$ in $\pi_0^{-1}(p)$ equals the
codimension of the subvariety
\begin{equation}
\beta(\mathcal{Y}) \subset \prod_{j=1}^d
\PP\lt( \OO_{\PP^n}(d-j)\otimes_\OO \text{Sym}^j Q^\vee \rt)|_p.
\end{equation} 
By construction, $\beta(\mathcal{Y})$ is the locus parametrizing
sequences of hypersurfaces in $\PP Q|_p$, $(X_1,\dots,X_d)$, of
degrees $1,\dots, d$ respectively such that the intersection
\begin{equation}
X_{(1,\dots,d)} := X_1\cap \dots \cap X_d
\end{equation}
has dimension greater than $n-d-1$.  So we have reduced
theorem~\ref{thm-thm1} to the following theorem:

\begin{thm}\label{thm-thm2}  Let $Q$ be a vector space over $k$ of
  dimension $n$ and let $d$ be an integer such that $1\leq d \leq
  n-1$.  Let $\PP_d$ denote the scheme $\prod_{j=1}^d\PP \text{Sym}^j
  Q^\vee$, which parametrizes $d$-tuples $\lt(X_1,\dots,X_d\rt)$ of
  hypersurfaces $X_i\in \PP Q$ of degree $i$.  Denote by $D_d$ the
  closed subscheme of $\PP_d$ which parametrizes sequences
  $\lt(X_1,\dots,X_d\rt)$ such that
\begin{equation}
\text{dim}\lt(X_{(1,\dots,d)}\rt)> n-d-1.
\end{equation}
The codimension of $D_d$ in
$\PP_d$ is greater than $n-1$.
\end{thm}

\begin{proof}
  We will prove this by induction on $d$.  
Consider first the case $d=1$.
Since $D_1=\emptyset$ and
  the dimension of $\PP_d=\PP Q^\vee$ is $n-1$, the result is true for
  $d=1$.

\

Let $U_d$ denote the open subscheme of $\PP_d$ which is the complement
of $D_d$.  Then for $1\leq d \leq n-2$, $U_{d+1}$ is contained in
$U_d\times \PP \text{Sym}^{d+1} Q^\vee$.  To see this, note that if
$X_{(1,\dots,d)}$ has dimension larger than $n-d-1$, then
$X_{1,\dots,d+1}$ is nonempty and has dimension greater than $n-d-2$:
it is nonempty since $X_{d+1}$ is ample, it has dimension larger than
$n-d-1$ by the Hauptidealsatz.  So we see that the codimension of
$D_{d+1}$ in $\PP_{d+1}$ is the minimum of the codimension of $D_d$ in
$\PP_d$ and the codimension of $D_{d+1}\cap (U_d\times \PP
\text{Sym}^{d+1} Q^\vee)$ in $U_d\times \PP \text{Sym}^{d+1} Q^\vee$.
So by induction we are reduced to showing that the codimension of
$D_{d+1}\cap (U_d\times \PP \text{Sym}^{d+1} Q^\vee)$ in $U_d\times
\PP \text{Sym}^{d+1} Q^\vee$ is larger than $n-1$.

\

Now suppose that $(X_1,\dots,X_d,X_{d+1})$ is a point in $D_{d+1}\cap
(U_d\times \PP \text{Sym}^{d+1} Q^\vee)$.  By assumption every
irreducible component of $X_{(1,\dots,d)}$ has dimension $n-d-1$.
Since also $X_{(1,\dots,d+1)}$ has dimension $n-d-1$, we conclude that
there is an irreducible component $C\subset X_{(1,\dots,d)}$ such that
$C\subset X_{d+1}$.  If $X_{(1,\dots,d)}=C_1\cup \dots \cup C_r$ is
the irreducible decomposition, then the fiber of $D_{d+1}\cap
(U_d\times \PP\text{Sym}^{d+1} Q^\vee)$ over $(X_1,\dots,X_d)$ (which
we consider as a subscheme of $\PP\text{Sym}^{d+1}Q^\vee$) is just the
union of $i=1,\dots,r$ of the set $B_i\subset \PP\text{Sym}^{d+1}
Q^\vee$ parametrizing hypersurfaces $X_{d+1}$ such that $C_i\subset
X_{d+1}$.  We are reduced to showing that the codimension of each
$B_i$ in $\PP\text{Sym}^{d+1} Q^\vee$ is greater than $n-1$.  We prove
this in a lemma:

\

\begin{lem}\label{lem-lem2}  Let $Y\subset \PP Q$ be an irreducible
  subscheme such that $\text{dim}Y=n-d-1$.  Let $B(Y)\subset
  \PP\text{Sym}^{d+1} Q^\vee$ be the locus of hypersurfaces $X_{d+1}$
  such that $Y\subset X_{d+1}$.  The codimension of $B(Y)$ is greater
  than $n-1$.
\end{lem}

\begin{proof}  Let $\Lambda\subset \PP Q$ be a $(d-1)$-plane which is disjoint
  from $Y$.  Choose coordinates on $\PP Q$, $(x_0,\dots,x_{n-1})$ with
  respect to which $\Lambda = Z(x_d,\dots,x_{n-1})$.  Let ${\mathbb
    G}_m$ denote the multiplicative group $\SP{\CC[t,t^{-1}]}$.  Let
  $\mu:{\mathbb G}_m\times \PP Q\rightarrow \PP Q$ be the group action
  given by
\begin{equation} 
\mu(t, (x_0, \dots, x_{d-1}, x_d, \dots, x_{n-1})) = (t^{-1}x_0,
\dots, t^{-1}x_{d-1}, tx_d, \dots, tx_{n-1}).
\end{equation}
Since the Hilbert scheme of $\PP Q$ is proper, the valuative criterion
implies that the closed subscheme
\begin{equation}
\mu^{-1}(Y)\subset {\mathbb G}_m\times \PP Q
\end{equation} 
which is flat over ${\mathbb G}_m$, extends over $0$ to yield a closed
subscheme
\begin{equation} 
\mathcal{Y}\subset \AAA^1\times \PP Q
\end{equation} 
which is flat over $\AAA^1$.  It is easy to see that the fiber of
$\mathcal{Y}$ over $0$ is a scheme whose reduced scheme is just
\begin{equation}
Z(x_0,\dots,x_{d-1})\subset \PP Q.
\end{equation}

\

Now we can form the family
\begin{equation}
\mathcal{B}\subset \AAA^1\times \PP\text{Sym}^{d+1} Q^\vee,
\mathcal{B}_t=B(\mathcal{Y}_t).
\end{equation} 
Over $\mathbb{G}_m$ the fibers of $\mathcal{B}$ are isomorphic.  It
follows by upper semicontinuity that for $t\neq 0$ we have
$\text{dim}(\mathcal{B}_t)\leq \text{dim}(\mathcal{B}_0)$.  And of
course we have
\begin{equation}
\mathcal{B}_0=B(\mathcal{Y}_0)\subset
B \lt(Z(x_0,\dots,x_{d-1}) \rt).
\end{equation}  
So we are reduced to proving the lemma for the special case
$Y=Z(x_0,\dots,x_{d-1})$.  The set B of hypersurfaces $X_{d+1}\subset
\PP Q$ which contain $Z(x_0,\dots,x_{d-1})$ is just the
projectivization of the kernel of the surjective linear map
$$H^0(\PP Q,\OO_{\PP Q}(d+1))\rightarrow H^0(Y,\OO_Y(d+1)).$$  
So the codimension of B in $\PP Q$ equals 
$$\text{dim}_{\CC} H^0(Y,\OO_Y(d+1))=\binom{n}{d+1}.$$  
For $d+1\leq n-1$ (which is one of
our hypotheses) we see that $\binom{n}{d+1}\geq n > n-1$.  We conclude
that the codimension of $B(Y)$ in $\PP\text{Sym}^{d+1}Q^\vee$ is
greater than $n-1$.  This proves the lemma.
\end{proof}

By the above lemma, we conclude that the codimension of each $B_i$ in
$\PP \text{Sym}^{d+1} Q^\vee$ is greater than $n-1$.  So we have
proved theorem~\ref{thm-thm2}
\end{proof}

Since we had reduced theorem~\ref{thm-thm1} to theorem~\ref{thm-thm2},
we have proved theorem~\ref{thm-thm1}.

\

While we are discussing results about lines on hypersurfaces, let us
mention two other results about lines on hypersurfaces.

\begin{lem} ~\cite[exercise V.4.4.2]{K} \label{lem-dfmlines}  For
  general $X$ and a general line $L\subset X$, the normal bundle
  $N_{L/X}$ is of the form $\OO_L^{\oplus d-1}\oplus\OO_L(1)^{\oplus
    n-1-d}$.
\end{lem}

\begin{thm} ~\cite[theorem V.4.3.2]{K} \label{lem-smlines}  For general
  X, the Fano scheme $F_1(X)$ is smooth.  Therefore $F_{0,1}(X)$ is
  smooth.  By generic smoothness, the general fiber of
  $F_{0,1}(X)\rightarrow X$ is smooth.
\end{thm}

\section{Stable $A$-graphs and Stable Maps}\label{sec-stab}

We follow the notation from ~\cite{BM} regarding stable A-graphs.
However, we shall only need to use genus 0 trees.

\subsection{Graphs and Trees}~\label{subsec-graph}

\begin{defn}~\label{def-graph}
  A \emph{graph} $\tau$ is a $4$-tuple
  $(F_\tau,W_\tau,j_\tau,\partial_\tau)$ defined as follows:
\begin{enumerate}
\item $F_\tau$ is a finite set called the set of \emph{flags}
\item $W_\tau$ is a finite set called the set of \emph{vertices}
\item $j_\tau:F_\tau\rightarrow F_\tau$ is an involution
\item $\partial_\tau:F_\tau\rightarrow W_\tau$ is a map called the
\emph{evaluation} map.
\end{enumerate}
In addition we have the auxiliary definitions
\begin{enumerate}
\item the set of \emph{tails} $S_\tau\subset F_\tau$ is the set of
fixed points of $j_\tau$
\item the set of \emph{edges} $E_\tau$ is the quotient of
$F_\tau\setminus S_\tau$ by $j_\tau$
\item for a vertex $v\in W_\tau$, the \emph{valence} of $v$ is defined
to be $\text{val}(v)=\#(\partial^{-1}(v))$.
\end{enumerate}
We shall often write $\text{Flag}(\tau)$ in place of $F_\tau$,
$\text{Vertex}(\tau)$ in place of $W_\tau$, $\text{Tail}(\tau)$ in
place of $S_\tau$, $\text{Edge}(\tau)$ in place of $E_\tau$, and
$\overline{f}$ in place of $j_\tau(f)$.
\end{defn}

We can associate to a graph its \emph{geometric realization} $|\tau|$
which is a simplicial complex defined as follows.  The set of
0-simplices of $|\tau|$ is
\begin{equation}
|\tau|^0 = \text{Vertex}(\tau) \sqcup
\text{Tail}(\tau).
\end{equation}  
The set of 1-simplices of $|\tau|$ is
\begin{equation}
|\tau|^1 = \text{Edge}(\tau)\sqcup \text{Tail}(\tau).
\end{equation}
If $[0,1]$ is a 1-simplex associated to an edge
$\lt\{f,\overline{f}\rt\}$, the point 0 is glued to the 0-simplex
$\partial f$, and the point 1 is glued to the 0-simplex $\partial
\overline{f}$.  If $[0,1]$ is the 1-simplex associated to a tail $f$,
the point 0 is glued to the 0-simplex $\partial f$, and the point 1 is
glued to the 0-simplex $f$.

\begin{defn}~\label{def-tree}
  A \emph{tree} is a connected graph such that $H_1(|\tau|,\ZZ) = 0$,
  i.e. a graph which contains no closed loops.
\end{defn}

One important tree is the \emph{empty tree} $\lambda_\emptyset$, i.e.
the tree such that $\text{Vertex}(\lambda_\emptyset)=\emptyset$.  For
each nonnegative integer $r$ define $\lambda_r$ to be the tree with
one vertex, $\text{Vertex}(\lambda_r)=\{v\}$, and with $r$ flags (all
of which are tails), $\text{Tail}(\lambda_r)=\{f_1,\dots,f_r\}$.
Also, for each pair of nonnegative integers $(r_1,r_2)$, define
$\lambda_{r_1,r_2}$ to be the connected tree with two vertices
$v_1,v_2$, with $r_1$ tails attached to $v_1$ and with $r_2$ tails
attached to $v_2$.

\begin{defn}\label{def-Agraph} 
  An \emph{$A$-graph} is a pair $(\tau,\beta_\tau)$ where $\tau$ is a
  tree and
\begin{equation}
\beta:\text{Vertex}(\tau)\rightarrow \ZZ_{\geq 0}
\end{equation}
is a map called the \emph{$A$-structure}.  We shall often abbreviate
$(\tau,\beta_\tau)$ by just writing $\tau$.  We say that an $A$-graph
$\tau$ is \emph{stable} if for each vertex $v\in\text{Vertex}(\tau)$
such that $\beta_\tau(v)=0$, there are at least 3 distinct flags
$f\in\text{Flag}(\tau)$ such that $\partial f = v$ (i.e. the valence
of v is at least 3).
\end{defn}

One important $A$-graph is the empty graph, $\tau_\emptyset$.  This is
the unique $A$-graph whose underlying graph is $\lambda_\emptyset$.
For each pair of nonnegative integers $r$ and $e$, define $\tau_r(e)$
to be the unique $A$-graph whose underlying graph is $\lambda_r$ and
such that $\beta(v)=e$.  Obviously $\tau_r(e)$ is stable iff either
$r\geq 3$ or $e> 0$.  For each pair of pairs $(r_1,r_2)$ and
$(e_1,e_2)$ where $r_1,r_2,e_1$ and $e_2$ are nonnegative integers,
define $\tau_{r_1,r_2}(e_1,e_2)$ to be the unique $A$-graph whose
underlying graph is $\lambda_{r_1,r_2}$, such that $\beta(v_1)=e_1$
and such that $\beta(v_2)=e_2$.

\

There is a category whose objects are the stable $A$-graphs.  The
morphisms in this category are each a composition of two basic types
of morphisms: \emph{contractions} and \emph{combinatorial morphisms},
c.f. ~\cite{BM} for the precise definitions.  Essentially a
contraction of $A$-graphs $\phi:\tau\rightarrow \sigma$ is a map from
the vertices of $\tau$ onto the vertices of $\sigma$ which maps
adjacent vertices to adjacent vertices (here two vertices are adjacent
if they are equal or if they are connected by an edge).  And a
combinatorial morphism $\tau \hookleftarrow \sigma$ is the inclusion
of a subgraph $\sigma$ into a graph $\tau$.  The functor which
associates to a stable $A$-graph the corresponding \emph{Behrend-Manin
  stack} is covariant for contractions.  But it is contravariant for
combinatorial morphisms.  Therefore we think of a combinatorial
morphism $\tau \hookleftarrow \sigma$ as a morphism from $\tau$ to
$\sigma$ (which explains our terminology $\tau \hookleftarrow \sigma$
for combinatorial morphisms).

\

Particularly important are morphisms of graphs which remove tails.
For each stable $A$-graph $\tau$ we define $r_{>0}(\tau)$ to be the
stable $A$-graph obtained by removing every tail
$f\in\text{Tail}(\tau)$ such that $\beta(\partial f) > 0$.  We define
$\tau \hookleftarrow r_{>0}(\tau)$ to be the canonical combinatorial
morphism.  For each stable $A$-graph $\tau$ we define $r_0(\tau)$ to
be the stabilization of the $A$-graph obtained by removing all tails
$f\in\text{Tail}(\tau)$ such that $\beta(\partial f)=0$.  Technically
the canonical morphism of graphs from $\tau$ to $r_0(\tau)$ consists
of both a combinatorial morphism and a contraction.  But we shall
denote it by $\tau \hookleftarrow r_0(\tau)$ just as if it were a
combinatorial morphism.  Finally, we define $r(\tau) :=
r_{>0}(r_0(\tau))=r_0(r_{>0}(\tau))$.

\

There are numerical invariants associated to an $A$-graph.

\begin{defn}~\label{def-invts}
  Given an $A$-graph $\tau$, define
\begin{equation}
\beta(\tau) = \sum_{v\in \text{Vertex}(\tau)} \beta(v).
\end{equation}
If $(X,L)$ is a polarized variety such that $K_X
\stackrel{\text{num}}{=} m L$ for some integer $m$, define the
\emph{expected dimension}
\begin{equation}
\text{dim}(X,\tau) = -m\beta(\tau) +
\#\text{Tail}(\tau) - \#\text{Edge}(\tau) + (\text{dim}(X)-3).
\end{equation}
\end{defn}

\

\subsection{Prestable Curves and Dual Graphs}~\label{subsec-dual}

\begin{defn}~\label{def-pcurve}  
  A \emph{prestable curve} with $r$ marked points
  $(C,(x_1,\dots,x_r))$ is a pair where $C$ is a complete, reduced, at
  worst nodal curve and $x_i\in C, i=1,\dots,r$ are distinct,
  nonsingular points of $C$.
\end{defn}

Suppose that $(C,(x_1,\dots,x_r))$ is a connected, prestable curve
whose arithmetic genus is $0$.  One associates to
$(C,(x_1,\dots,x_n))$ a \emph{dual graph}, $\Delta$: a tree whose
vertices $\{v_1,v_2,\dots\}$ correspond to the irreducible components
$\{C_1, C_2,\dots \}$ of $C$, whose edges $\{ \{f_1,\overline{f_1}\},
\{f_2,\overline{f_2}\},\dots \}$ correspond to the nodes
$\{q_1,q_2,\dots \}$ of $C$, and whose tails $\{ g_1,\dots, g_r\}$
correspond to the marked points $\{p_1,\dots,p_r\}$ of $C$.

\begin{defn}~\label{def-pmap} 
  Let $(X,L)$ be a polarized variety.  A \emph{prestable map} is a
  pair
\begin{equation}
((C,(x_1,\dots,x_n)),C\xrightarrow{h} X)
\end{equation}
where $(C,(x_1,\dots,x_n))$ is a prestable curve, and where
$C\xrightarrow{h} X$ is a morphism of $\CC$-schemes.
\end{defn}

Just as one associates to a connected prestable curve
$(C,x_1,\dots,x_n)$ of arithmetic genus $0$ a tree $\Delta(C,x)$, one
can associate an $A$-graph to a prestable map $\lt( (C,x_1,\dots,x_n),
C\xrightarrow{h} X\rt)$ from a connected prestable curve of arithmetic
genus $0$.  The underlying tree of $\Delta(C,x,h)$ is simply
$\Delta(C,x)$.  And, given a component $C_i$ of $C$ with corresponding
vertex $v_i\in \text{Vertex}(\Delta(C,x))$, one defines
\begin{equation}
\beta(v_i)=\int_{C_i} h_i^*(c_1(L)).
\end{equation}
The $A$-graph $\Delta(C,x,h)$ is a \emph{stable} $A$-graph iff
$(C,x,h)$ is a \emph{stable} map.

\subsection{Behrend-Manin stacks}~\label{subsec-BM}
  
We refer the reader to ~\cite{BM} for the definition of the stacks
$\Kbm{X}{\tau}$.  These are proper Deligne-Mumford stacks which
parametrize stable maps along with some extra data.  We shall
sometimes deal with these stacks, but more often we shall deal with
the open substack $\Kbmo{X}{\tau}\subset \Kbm{X}{\tau}$ of
\emph{strict maps} which we now define.

\begin{defn}~\label{defn-strict}
  Let $X$ be a variety, $L$ a line bundle on $X$, and let $\tau$ be a
  stable A-graph.  A \emph{strict $\tau$-map} is a datum
\begin{equation}
((C_v),(h_v:C_v\rightarrow X),(q_f))
\end{equation}
defined as follows:
\begin{enumerate}
\item $(C_v)$ is a set parametrized by $v\in\text{Vertex}(\tau)$ of
rational curves, i.e. each $C_v\cong \PP^1$
\item $(h_v:C_v\rightarrow X)$ is a set parametrized by
$v\in\text{Vertex}(\tau)$ of morphisms of $\CC$-schemes,
\item $(q_f)$ is a set parametrized by $f\in\text{Flag}(\tau)$ of
closed points $q_f\in C_{\partial f}$
\end{enumerate}
and satisfying the following conditions
\begin{enumerate}
\item for $v\in\text{Vertex}(\tau)$, the degree of $h_v^*(L)$ as a
line bundle on $C_v$ is $\beta_\tau(v)$,
\item for $f_1,f_2\in\text{Flag}(\tau)$ distinct flags with $\partial
f_1=\partial f_2$, $q_{f_1}\neq q_{f_2}$,
\item for $f\in\text{Flag}(\tau)$, we have $h_{\partial
f}(q_f)=h_{\partial \overline{f}}(q_{\overline{f}})$.
\end{enumerate}
\end{defn}

\textbf{Convention:} For the empty graph, $\tau_\emptyset$, we define a
strict $\tau_\emptyset$-map to simply be a point in $X$.  Thus the set
of strict $\tau_\emptyset$-maps is simply $X$.

\begin{defn}~\label{defn-strictfam}
  If $T$ is a $\CC$-scheme, then a \emph{family of strict $\tau$-maps}
  over $T$ is a datum
\begin{equation}
((\pi_v:\mathcal{C}_v\rightarrow T),
(h_v:\mathcal{C}_v\rightarrow X), (q_f:T\rightarrow
\mathcal{C}_{\partial f}))
\end{equation}
defined as follows:
\begin{enumerate}
\item $(\pi_v:\mathcal{C}_v\rightarrow T)$ is a set parametrized by
$v\in\text{Vertex}(\tau)$ of smooth, proper morphisms whose geometric
fibers are rational curves
\item $(h_v:\mathcal{C}_v\rightarrow X)$ is a set parametrized by
$v\in\text{Vertex}(\tau)$ of morphisms of $\CC$-schemes
\item $(q_f:T\rightarrow \mathcal{C}_{\partial f})$ is a set parametrized by
$f\in\text{Flag}(\tau)$ of morphisms of schemes such that
$\pi_{\partial f}\circ q_f = id_T$
\end{enumerate}
and satisfying the following conditions
\begin{enumerate}
\item for $v\in\text{Vertex}(\tau)$, the degree of $h_v^*(L)$ on each
geometric fiber of $\mathcal{C}_v\rightarrow S$ is $\beta_\tau(v)$
\item for $f_1,f_2\in\text{Flag}(\tau)$ distinct flags with $\partial
f_1 = \partial f_2$, $q_{f_1}$ and $q_{f_2}$ are disjoint sections
\item for $f\in\text{Flag}(\tau)$, we have $h_{\partial f}\circ q_f =
h_{\partial \overline{f}}\circ q_{\overline{f}}$.
\end{enumerate}
\end{defn}

\textbf{Convention:} For the empty graph, $\tau_\emptyset$ we define a
family of strict $\tau_\emptyset$-maps over $T$ to be a morphism
$h:T\rightarrow X$.

\

Suppose given two families of strict $\tau$-maps over $S$, say
\begin{eqnarray}
\eta= ((\pi_v:\mathcal{C}_v\rightarrow T),
(h_v:\mathcal{C}_v\rightarrow X), (q_f:T\rightarrow
\mathcal{C}_{\partial f})), \\
\zeta=
((\pi'_v:\mathcal{C}'_v\rightarrow T),
(h'_v:\mathcal{C}'_v\rightarrow X), (q'_f:T\rightarrow
\mathcal{C}'_{\partial f})).
\end{eqnarray}

\begin{defn}~\label{defn-strictmor}
  A \emph{morphism} of families of strict $\tau$-maps over $S$,
  $\phi:\eta\rightarrow \zeta$, is a collection of isomorphisms of
  $S$-schemes:
\begin{equation}
\phi =
(\phi_v:\mathcal{C}_v\rightarrow \mathcal{C}'_v)
\end{equation}
indexed by
$v\in\text{Vertex}(\tau)$ and satisfying
\begin{enumerate}
\item for $v\in\text{Vertex}(\tau)$, $h'_v\circ \phi_v = h_v$
\item for $f\in\text{Flag}(\tau)$, $\phi_{\partial f}\circ q_f = q'_f$.
\end{enumerate}
\end{defn}

One defines composition of morphisms in the obvious way.  Notice that
every morphism is an isomorphism.  Thus the category of families of
strict $\tau$-maps over $S$ is a groupoid.  Given a morphism
$S'\xrightarrow{u} S$ and a family $\eta$ of strict $\tau$-maps over
$S$, one has the usual pullback $u^*(\eta)$ which is a family of
strict $\tau$-maps over $S'$.  In this way we have the notion of a
functor from the category of $\CC$-schemes to the category of
groupoids which associates to each $S$ the groupoid of families of
strict $\tau$-maps over $S$.  We denote this functor by
$\Kbmo{X}{\tau}$.

\

In every case it is easy to see that $\Kbmo{X}{\tau}$ is a stack in
groupoids over $\CC$.  In many cases this is even a Deligne-Mumford
stack:

\begin{thm}~\label{thm-KBMO}  If $X$ is projective and $L$ is ample,
  the functor $\Kbmo{X}{\tau}$ is a Deligne-Mumford stack which is
  separated and finite type over $\CC$.
\end{thm}

\begin{proof}
  There is a 1-morphism $\Kbmo{X}{\tau}\rightarrow \Kbm{X}{\tau}$
  where $\Kbm{X}{\tau}$ is the refined functor defined in ~\cite{BM}.
  In ~\cite{BM} it is proved that $\Kbm{X}{\tau}$ is a proper
  Deligne-Mumford stack over $\CC$.  And it is clear that
  $\Kbmo{X}{\tau}\rightarrow \Kbm{X}{\tau}$ is a representable
  morphism which is an open immersion.  Thus $\Kbmo{X}{\tau}$ is a
  Deligne-Mumford stack which is separated and finite type over $\CC$.
\end{proof}

\subsection{Properties and Related Constructions}~\label{subsec-props}

Notice that with our notation $\Kbm{X}{\tau_r(e)}$ is the moduli stack
of Kontsevich stable maps $\Kgnb{0}{r}{X}{e}$, and
$\Kbmo{X}{\tau_r(e)}$ simply parametrizes those stable maps such that
the domain curve is irreducible.

\

\begin{defn}~\label{defn-evmor}
  Suppose that $\tau$ is a stable $A$-graph and
  $f\in\text{Flag}(\tau)$.  Then there is a $1$-morphism
\begin{equation}
\text{ev}_f:\Kbm{X}{\tau}\rightarrow X
\end{equation}
defined by sending a family of $\tau$-maps, $\eta$ (with notation as
above), to the morphism $h_{\partial f}\circ q_f$.
\end{defn}

If $\alpha=(\alpha_F,\alpha_V):\sigma\rightarrow\tau$ is a
combinatorial morphism of 
graphs, $\tau\hookleftarrow \sigma$, there is an associated
$1$-morphism
\begin{equation}
        \Kbm{X}{\alpha}:\Kbm{X}{\tau}\rightarrow\Kbm{X}{\sigma}.
\end{equation}
If $\alpha$ is the inclusion of $\sigma$ as a subgraph of $\tau$, then
$\Kbm{X}{\alpha}$ is the \emph{forgetful morphism} which ``remembers'' only
those components of $\tau$-maps whose vertex is contained in $\sigma$.
The reader is referred to~\cite[theorem 3.6]{BM} for the precise
definition.  We will refer to the restriction of $\Kbm{X}{\alpha}$ to
$\Kbmo{X}{\tau}$ as $\Kbmo{X}{\alpha}$.

\

If $\phi=(\phi_W,\phi^F):\tau\rightarrow \tau'$ is a contraction of
stable $A$-graphs, there is a corresponding 1-morphism of proper
Deligne-Mumford stacks
\begin{equation}
\Kbm{X}{\phi}:\Kbm{X}{\tau}\rightarrow \Kbm{X}{\tau'}.
\end{equation}
This morphism ``forgets'' the labeling of some the individual
components of the domain curve.  The reader is referred
to~\cite[theorem 3.6]{BM} for the precise definition.  We will denote
by $\Kbmo{X}{\phi}$ the restriction of this 1-morphism to the open
substack $\Kbmo{X}{\tau}$ of $\Kbm{X}{\tau}$.

\

One important case to understand is when $\beta(\tau)=0$.  We have
already defined $\Kbmo{X}{\tau_\emptyset}=\Kbm{X}{\tau_\emptyset}=X$
where $\tau_\emptyset$ is the empty graph.  For any stable $A$-graph
$\tau$ such that $\beta(\tau)=0$ and such that
$\#\text{Tail}(\tau)=r$, we have $\Kbmo{X}{\tau}= X \times
\Kbmo{*}{\tau}$ where $\Kbmo{*}{\tau}\subset \Mgnb{0}{r}$ is the
obvious substack.

\

Consider the case that $\phi:\tau\rightarrow \tau'$ is a contraction
of stable A-graphs such that $\beta(\tau)=\beta(\tau')=0$.  Then
$\Kbmo{X}{\tau}\rightarrow \Kbm{X}{\tau'}$ is simply the product of
$\text{id}_X:X\rightarrow X$ with the 1-morphism
$\Kbmo{*}{\phi}:\Kbmo{*}{\tau}\rightarrow \Kbm{*}{\tau'}$.  In
particular, using the notation of ~\cite{BM}, consider the case that
$\phi$ is an \emph{isogeny}, i.e. $\phi$ is the morphism which removes
some subset of the set of tails from $\tau$ and then stabilizes the
resulting (possibly unstable) graph.

\begin{lem}~\label{lem-isog}  Let $\tau,\tau'$ be stable $A$-graphs such
  that $\beta(\tau)=\beta(\tau')=0$ and let $\phi:\tau\rightarrow
  \tau'$ be an isogeny.  Then
  $\Kbmo{X}{\phi}:\Kbmo{X}{\tau}\rightarrow \Kbm{X}{\tau'}$ is smooth
  of relative dimension $\text{dim}(X,\tau)-\text{dim}(X,\tau')$ with
  geometrically connected fibers.
\end{lem}

\begin{proof}  Of course it is equivalent to prove that
\begin{equation}
\Kbmo{*}{\phi}:\Kbmo{*}{\tau}\rightarrow \Kbm{*}{\tau'}
\end{equation}
is smooth of relative 
dimension $\text{dim}(X,\tau)-\text{dim}(X,\tau')$ with geometrically
connected fibers.  Now it follows by proposition 7.4 of
~\cite{BM} that $\Kbmo{*}{\tau}$ and $\Kbm{*}{\tau'}$ have the
expected dimension.  Thus all we really need to show is that
$\Kbmo{*}{\tau}\rightarrow \Kbm{*}{\tau'}$ is smooth with geometrically
irreducible fibers.  Moreover, since every isogeny is a composition of
morphisms obtained by stably removing one tail, we may suppose that
$\phi:\tau\rightarrow \tau'$ corresponds to stably removing one tail
$f\in\text{Tail}(\tau)$.  

\

There are two cases.  Suppose first of all that when we remove $f$
from $\tau$, the resulting graph is unstable.  But then $\partial f =
v$ is a vertex with valence 3.  Since a rational curve with 3 marked
points has no moduli, we conclude that $\Kbmo{*}{\tau}\rightarrow
\Kbm{*}{\tau'}$ is an open immersion.

\

The second case is that when we remove $f$ from $\tau$, the resulting
graph is stable, i.e. the resulting graph is just $\tau'$.  But then
if $v=\phi(\partial f)$, we conclude that
$\phi:\Kbmo{*}{\tau}\rightarrow \Kbm{*}{\tau'}$ is simply an open
subset of the universal curve over $\Kbm{*}{\tau'}$ corresponding to
the vertex $v$.  In both cases we conclude that
$\Kbm{*}{\tau}\rightarrow \Kbm{*}{\tau'}$ is smooth with geometrically
connected fibers.
\end{proof}

For each stable $A$-graph $\tau$ define $e=\beta(\tau)$ and define
$r=\#\text{Tail}(\tau)$.  Then there is a contraction
$\phi:\tau\rightarrow \tau_r(e)$ which is unique up to a labeling of
the tails of $\tau$.

\begin{defn}~\label{defn-cancon}
  The contraction $\phi$ above is the \emph{canonical contraction}.
  The corresponding $1$-morphism (well-defined up to relabeling the
  tails)
\begin{equation}
\Kbmo{X}{\phi}:\Kbmo{X}{\tau}\rightarrow\Kgnb{0}{r}{X}{e}
\end{equation}
will be referred to as the \emph{canonical contraction morphism}.
\end{defn}

Notice that the image of $\Kbmo{X}{\phi}$ as a subset of the set
$|\Kgnb{0}{r}{X}{e}|$ is well-defined.

\begin{prop}~\label{prop-locclo}
  Let $\phi:\tau\rightarrow \tau'$ be a contraction of stable
  $A$-graphs.  The image of the 1-morphism $\Kbmo{X}{\phi}$ is a
  locally closed subset of the topological space $|\Kbm{X}{\tau}|$.
\end{prop}

\begin{proof}  For notation's sake let's denote the continuous map
  of topological spaces
\begin{equation}
|\Kbm{X}{\phi}|:|\Kbm{X}{\tau}|\rightarrow
|\Kbm{X}{\tau'}|
\end{equation}
by $ f: M\rightarrow M'$ and let's denote the
open substack $\Kbmo{X}{\tau}$ of $\Kbm{X}{\tau}$ by $M^o$.  Then
$f:M\rightarrow M'$ is a closed map.  And it is easy to see that
$f^{-1}(f(M^o)) = M^o$.  Therefore $f(M^o)= f(M) - f(M - M^o)$ is
a difference of closed sets and so is locally closed.
\end{proof}

\

We now fix $n$ and $\alpha$ and consider the set $S$ of all images
\begin{equation}
\lt\{\text{Image}(\Kbmo{X}{\phi})\rt\}
\end{equation}
as $\phi$ ranges
over all contractions of stable $A$-graphs to $\tau_n(\alpha)$.
The set of isomorphism classes of such contractions is clearly
finite.  The previous lemma shows that $S$ forms a \emph{locally
closed decomposition} of the topological space
$|\Kgnb{0}{n}{X}{\alpha}|$, i.e. a partition of
$|\Kgnb{0}{n}{X}{\alpha}|$ into locally closed subsets.  This
partition is what we call the \emph{Behrend-Manin decomposition}.

\section{Flatness and Dimension Results}~\label{sec-flat}

In this section we consider the dimensions of the stacks
$\Kbmo{X}{\tau}$ and the evaluation morphisms.  The main property we
are interested in is the following:

\begin{defn}~\label{defn-prop1}
  Given a stable $A$-graph $\tau$, we say that $\cD(X,\tau)$ holds if
  the dimension of every irreducible component of $\Kbmo{X}{\tau}$
  equals the expected dimension $\text{dim}(X,\tau)$.
\end{defn}

By deformation theory there is an a priori lower bound on the
dimension of any irreducible component of $\Kbmo{X}{\tau}$:

\begin{lem}~\label{lem-lbdim}
  Every irreducible component of $\Kbm{X}{\tau}$ has dimension at
  least $\text{dim}(X,\tau)$.  In particular, every irreducible
  component of $\Kbmo{X}{\tau}$ has dimension at least
  $\text{dim}(X,\tau)$.
\end{lem}

\begin{proof}
  This is a standard result. In the case that $\tau=\tau_r(e)$ it
  follows from ~\cite[section 5.2]{FP}.  In the general case the
  theorem follows from ~\cite{B}.  The theorem isn't actually stated
  in ~\cite{B}, so we show how it follows from results there.

\

Let $\mathfrak{M}(\tau)$ denote the stack of $\tau$-marked prestable
curves as in ~\cite{B}.  There is a forgetful $1$-morphism of
algebraic (Artin) stacks
\begin{equation}
        \Kbm{X}{\tau}\rightarrow \mathfrak{M}(\tau).
\end{equation}
By ~\cite[lemma 1]{B}, $\mathfrak{M}(\tau)$ is smooth of dimension
\begin{equation}
        \#\text{Tail}(\tau) - \#\text{Edge}(\tau) - 3.
\end{equation}
Let $\mathfrak{C}(\tau)\rightarrow \mathfrak{M}(\tau)$ be the
universal curve.  By ~\cite[proposition 4]{B}, $\Kbm{X}{\tau}$ is an
open substack of the relative morphism-scheme,
\begin{equation}
\text{Mor}_{\mathfrak{M}(\tau)}\lt(\mathfrak{C}(\tau),X\times
\mathfrak{M}(\tau)\rt).
\end{equation}
By ~\cite[theorem 2.17.1]{K}, it follows that every irreducible
component of $\Kbm{X}{\tau}$ has dimension at least
\begin{equation}
        h^0(C,f^*T_X)-h^1(C,f^*T_X)
        +\#\text{Tail}(\tau)-\#\text{Edge}(\tau). 
\end{equation}
By Riemann-Roch, we have $(h^0-h^1)(C,f^*T_X) = -K_X.f_*[C] +
\text{dim}(X)$.  We conclude that every irreducible component of
$\Kbm{X}{\tau}$ has dimension at least
\begin{equation}
        -K_X.f_*[C] + \text{dim}(X) - 3 + \#\text{Flag}(\tau) -
         \#\text{Edge}(\tau) = \text{dim}(X,\tau).
\end{equation}
\end{proof}

When $\text{Vertex}(\tau)$ has more than one element, we can try to
reduce $\cD(X,\tau)$ to $\cD(X,\tau_i)$ for some proper subgraphs
$\tau_i$ of $\tau$, thus giving an inductive proof that $\cD(X,\tau)$
holds.  To carry out such a proof, we need to know that the evaluation
morphisms have constant fiber dimension.  So we introduce the
following property:

\begin{defn}~\label{defn-prop2}
  Given a stable $A$-graph $\tau$ and a flag $f\in\text{Flag}(\tau)$,
  we say that $\cE(X,\tau,f)$ holds if
\begin{equation}
\text{ev}_f:\Kbmo{X}{\tau}\rightarrow X
\end{equation}
is dominant and has constant fiber dimension
$\text{dim}(X,\tau)-\text{dim}(X)$.  
\end{defn}

Notice that if there is any flag $f\in\text{Flag}(\tau)$ such that
$\cE(X,\tau,f)$ holds, then it follows that $\cD(X,\tau)$ holds.

\

In the case that $X\subset \PP^N$ is a complete intersection, then the
properties $\cD$ and $\cE$ are equivalent to stronger properties.

\begin{defn}~\label{defn-prop1'}
  Given a stable $A$-graph $\tau$, we say that $\lci(X,\tau)$ holds if
  $\Kbmo{X}{\tau}$ is a local complete intersection and if
  $\cD(X,\tau)$ holds.  Given a stable $A$-graph $\tau$ and a flag
  $f\in\text{Flag}(\tau)$, we say that $\fe(X,\tau,f)$ holds if
\begin{equation}
\text{ev}_f:\Kbmo{X}{\tau} \rightarrow X
\end{equation}
is flat of relative dimension $\text{dim}(X,\tau)-\text{dim}(X)$.
\end{defn}

\begin{lem}~\label{lem-propeq}
  If $X\subset \PP^N$ is a complete intersection, then $\cD(X,\tau)$
  holds iff $\lci(X,\tau)$ holds.  Also $\cE(X,\tau,f)$ holds iff
  $\fe(X,\tau,f)$ holds.  The same result holds with $\Kbmo{X}{\tau}$
  replaced by $\Kbm{X}{\tau}$.
\end{lem}

\begin{proof}  
  Suppose that $X$ is a complete intersection of $r = N-n$
  hypersurfaces of degrees $d_1,\dots,d_r$.  Consider
  $\Kbmo{\PP^N}{\tau}$ and denote the universal curve by
\begin{equation}
\pi:\mathcal{C}\rightarrow \Kbmo{\PP^N}{\tau}.
\end{equation}
Let $h:\mathcal{C}\rightarrow \PP^N$ denote the universal map.  Since
\begin{equation}
\OO_{\PP^N}({\mathbf d}) := \OO_{\PP^N}(d_1)\oplus \dots \oplus
\OO_{\PP^N}(d_r) 
\end{equation}
is generated by global sections, also $h^*\OO_{\PP^N}({\mathbf d})$ is
generated by global sections.  On a genus 0 tree, if $F$ is a sheaf
generated by global sections then $H^1(C,F)=0$, so
$R^1\pi_*\lt(h^*\OO_{\PP^N}({\mathbf d})\rt) = 0$.

\

Now by ~\cite[proposition 7.4]{BM}, $\Kbmo{\PP^N}{\tau}$ is smooth of
dimension
\begin{equation}
(N+1)\beta(\tau) + (N-3) +
\#\text{Flag}(\tau) - \#\text{Edge}(\tau).
\end{equation}
So by ~\cite[corollary II.2]{AVar} the pushforward
$E:=\pi_*\lt(h^*\OO_{\PP^N}({\mathbf d})\rt)$ is a locally free sheaf
of rank
\begin{equation}
 \sum_{i=1}^r \chi(C,h^*\OO_{\PP^N}(d_i)) =
\sum_{i=1}^r \lt(d_i\beta(\tau) + 1\rt) = \lt(\sum_i d_i\rt)\beta(\tau) + r.
\end{equation}
Now the defining equations of the hypersurfaces in $\PP^N$ give a
global section $\sigma$ of $E$, and $\Kbmo{X}{\tau}$ is precisely the
zero scheme of $\sigma$.  Finally notice that the expected
codimension, $\text{dim}(\PP^n,\tau) - \text{dim}(X,\tau)$, of
$\Kbmo{X}{\tau}$ in $\Kbmo{\PP^N}{\tau}$ is just
\begin{equation}
-K_{\PP^N}.f_*[C] +
K_X.f_*[C] + \text{dim}(\PP^N) - \text{dim}X = \lt(\sum_i
d_i\rt)\beta(\tau) + r.
\end{equation}
Thus, if $\Kbmo{X}{\tau}$, then it follows that $\Kbmo{X}{\tau}$ is a
local complete intersection.  So if $\cD(X,\tau)$ holds, then also
$\lci(X,\tau)$ holds.  The opposite inclusion is obvious.

\

Now suppose that $\cE(X,\tau,f)$ holds.  In particular $\cD(X,\tau)$
holds, so $\lci(X,\tau)$ holds.  But now by ~\cite[theorem 23.1]{Ma},
$\text{ev}_f$ is a dominant morphism from a Cohen-Macaulay scheme to a
smooth scheme with constant fiber dimension, therefore it is flat.  So
$\fe(X,\tau,f)$ holds.  The opposite inclusion is obvious.

\

The same proof works when we replace $\Kbmo{X}{\tau}$ by $\Kbm{X}{\tau}$.
\end{proof}

Consider the following diagram:

\begin{picture}(250,200)
\put(50,50){\makebox(200,100){\includegraphics{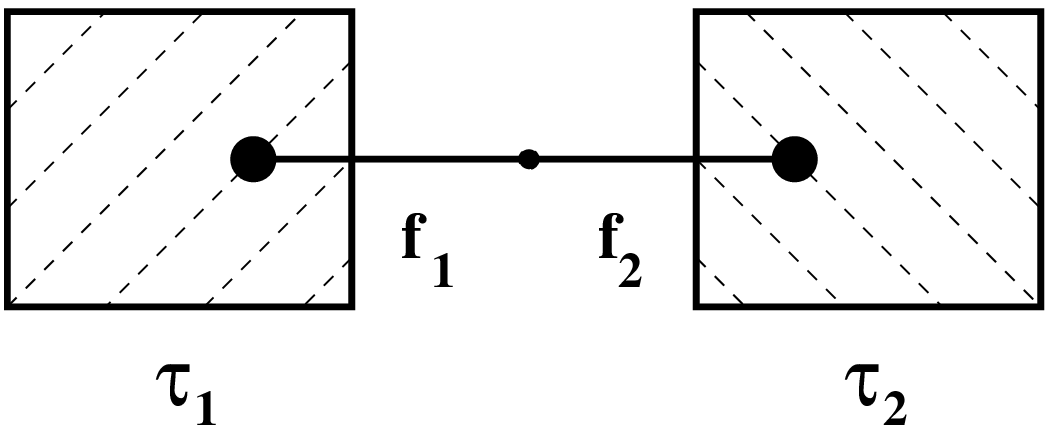}}}
\put(150,20){\makebox(0,0)[b]{Diagram 1}}
\end{picture}

Here $\tau$ is an $A$-graph which contains the two subgraphs 
\begin{equation}
a_1:\tau\hookleftarrow \tau_1,\ \  a_2:\tau\hookleftarrow \tau_2.
\end{equation}
The edge $\{f_1,f_2\}$ of $\tau$ is made up of the two tails
$f_1\in\text{Tail}(\tau_1), f_2\in\text{Tail}(\tau_2)$.  Let
$f\in\text{Flag}(\tau_2)$ be any flag (possibly $f=f_2$).

\begin{lem}~\label{lem-red1}
If $\fe(X,\tau_1,f_1)$ and $\fe(X,\tau_2,f)$
hold, then $\fe(X,\tau,f)$ holds.
\end{lem}

\begin{proof}
The combinatorial morphism $a_1,a_2$ give rise to a $1$-morphism 
\begin{equation}
\Kbmo{X}{a_1,a_2}:\Kbmo{X}{\tau}\rightarrow
\Kbmo{X}{\tau_1}\times_{\text{ev}_{f_1}, X, \text{ev}_{f_2}}
\Kbmo{X}{\tau_2}.
\end{equation}
It is clear (just by the definition of strict $\tau$-maps) that
$\Kbmo{X}{a_1,a_2}$ is an open immersion.  

\

Since $\text{ev}_{f_1}:\Kbmo{X}{\tau_1}\rightarrow X$ is flat of
relative dimension 
$\text{dim}(X,\tau_1) - \text{dim}(X)$, it follows that the projection
morphism 
\begin{equation}
\text{pr}_2: \Kbmo{X}{\tau_1}\times_{\text{ev}_{f_1},X, \text{ev}_{f_2}}
\Kbmo{X}{\tau_2} \rightarrow \Kbmo{X}{\tau_2}
\end{equation}
is flat of relative dimension $\text{dim}(X,\tau_1)-\text{dim}(X)$.
And $\text{ev}_{f}:\Kbmo{X}{\tau_2}\rightarrow X$ is flat of relative
dimension 
$\text{dim}(X,\tau_2) - \text{dim}(X)$.  Thus the composite morphism
\begin{equation}
\Kbmo{X}{\tau_1}\times_{\text{ev}_{f_1},X, \text{ev}_{f_2}}
\Kbmo{X}{\tau_2} \xrightarrow{\text{pr}_2} \Kbmo{X}{\tau_2}
\xrightarrow{\text{ev}_{f}} X
\end{equation}
is flat of relative dimension $\text{dim}(X,\tau_1) +
\text{dim}(X,\tau_2) - 2\text{dim}(X)$.
Of course $\text{ev}_f:\Kbmo{X}{\tau}\rightarrow X$ is simply the
restriction of the composite morphism, so it is flat of the same
relative dimension.
But notice that
\begin{equation}
\begin{array}{c}
 \text{dim}(X,\tau_1)+\text{dim}(X,\tau_2) - \text{dim}(X) = \\
 \lt(-K_X.\beta(\tau_1) + -K_X.\beta(\tau_2)\rt) +
2\lt(\text{dim}(X)-3\rt) + \\
\lt(\#\text{Flag}(\tau_1) +
\#\text{Flag}(\tau_2)\rt) + \lt(\#\text{Edge}(\tau_1) +
\#\text{Edge}(\tau_2)\rt) - \text{dim}(X)= \\
-K_X.\beta(\tau) + 2\lt(\text{dim}(X) - 3\rt) +
\lt(\#\text{Tail}(\tau) + 2\rt) - \\
\lt(\#\text{Edge}(\tau) - 1\rt)
-\text{dim}(X) = \text{dim}(X,\tau).
\end{array}
\end{equation}
From this it follows that  
\begin{equation}
\text{ev}_{f}:\Kbmo{X}{\tau}\rightarrow X
\end{equation}
is flat of constant fiber dimension
$\text{dim}(X,\tau) - \text{dim}(X)$.  Thus
$\fe(X,\tau,f)$ holds.  
\end{proof}

\begin{defn}~\label{defn-E}
Suppose that $\tau$ is a stable $A$-graph and define the
\emph{maximum component degree} of $\tau$ to be
\begin{equation}
E(\tau) = \sup_{v\in\text{Vertex}(\tau)} \beta(v).
\end{equation}
\end{defn}

\begin{prop}~\label{prop-red}
Suppose that $\tau$ is a stable $A$-graph with $E(\tau)=E$.  
If for each $e=0,\dots,E$ we have $\fe(X,\tau_1(e),f)$
holds, then for each flag $f\in\text{Flag}(\tau)$,
$\fe(X,\tau,f)$ holds.
\end{prop}

\begin{proof}
We prove this by induction on $\#\text{Vertex}(\tau)$.  Suppose $\tau$
has a single vertex.  Then $\tau=\tau_r(e)$ for some $r$ and $e$.  If
$e=0$, then $\Kbm{X}{\tau_r(0)}=\Kbm{*}{\tau_r(0)}\times X$, and the
evaluation morphism is just projection.  So $\text{ev}_f$ is flat of
relative dimension 
\begin{equation}
\text{dim}(\Kbm{*}{\tau_r(0)}) = \text{dim}(X,\tau_r(0)) -
\text{dim}(X),
\end{equation}
i.e.  $\fe(X,\tau_r(0),f)$ holds.  

\

Next consider the case $\tau=\tau_r(e)$ with 
$e>0$.  For any flag $f\in \text{Flag}(\tau_r(e))$ there is a
combinatorial morphism 
$a:\tau_{r}(e)\hookleftarrow \tau_1(e)$ which maps the unique flag 
$f_1\in\text{Flag}(\tau_1(e))$ to $f$.  The associated $1$-morphism
\begin{equation}
\Kbmo{X}{a}:\Kbmo{X}{\tau_r(e)}\rightarrow \Kbmo{X}{\tau_1(e)}
\end{equation}
is isomorphic to an open subset of the $(e-1)$-fold fiber product of
the universal curve over $\Kbmo{X}{\tau_1(e)}$.  Since the universal
curve is smooth of relative dimension $1$ over $\Kbmo{X}{\tau_1(e)}$,
we conclude that $\Kbmo{X}{a}$ is smooth of relative dimension $e-1$.  
The evaluation morphism
$\text{ev}_f:\Kbmo{X}{\tau_r(e)}\rightarrow X$ factors as the
composition
\begin{equation}
\begin{CD}
\Kbmo{X}{\tau_r(e)} @>\Kbmo{X}{a} >> \Kbmo{X}{\tau_1(e)} @>\text{ev}_{f_1}
>> X.
\end{CD}
\end{equation}
Therefore if $\fe(X,\tau_1(e),f_1)$ holds, then
$\text{ev}_f$ is flat of relative dimension
\begin{equation}
(e-1) + \lt(\text{dim}(X,\tau_1(e)) - \text{dim}(X)\rt) =
\text{dim}(X,\tau_r(e)) - \text{dim}(X),
\end{equation}
in other words $\fe(X,\tau_r(e),f)$ holds.
Thus the proposition is proved when $\tau$ has a single vertex.

\

Now suppose that $\tau$ has more than one vertex and suppose that for
all $e\leq E=E(\tau)$, $\fe(\tau_1(e),f)$ holds.  By
way of induction, assume that the
proposition is true for all graphs $\tau'$ such that
$\#\text{Vertex}(\tau') < \#\text{Vertex}(\tau)$.  
Let $\{f_1,f_2\}$
be any edge.  Define $\tau_1$ and $\tau_2$ to be the graphs obtained
by breaking the edge into two tails (see diagram 1).  
Let $f\in\text{Flag}(\tau)$ be a
flag, and without loss of generality suppose that
$f\in\text{Flag}(\tau_2)$.  Now $E(\tau_2) \leq E(\tau)$ and
$\#\text{Vertex}(\tau_2) < \#\text{Vertex}(\tau)$, so by the
induction assumption $\fe(X,\tau_2,f)$ holds.  Also
$E(\tau_1)\leq E(\tau)$ and $\#\text{Vertex}(\tau_1) <
\#\text{Vertex}(\tau)$, so by the induction assumption
$\fe(X,\tau_1,f_1)$ holds.  Then by lemma~\ref{lem-red1},
we conclude that $\fe(X,\tau,f)$ holds.  So the proposition
is proved by induction.
\end{proof}

\section{Specializations}

In the previous section we reduced the flatness and dimension results for
a general stable $A$-graph $\tau$ to flatness and dimension results
for the the stable $A$-graphs $\tau_1(e)$ with $0\leq e\leq
E(\tau).$
In this
section we will use specializations to reduce the flatness and
dimension results for all $\tau_1(e), e> 1$ to flatness and dimension
results for a finite number of cases $\tau_1(e), e=1,\dots,E(X)$ where
$E(X)$ is the \emph{threshold degree} of $X$.  We define a stable
$A$-graph $\sigma$ to be \emph{basic} if for each vertex
$v\in\text{Vertex}(\sigma)$, we have
$\beta(v)\leq E(X)$.  The
specializations we produce will show that every irreducible component
of $\Kbm{X}{\tau}$ contains a \emph{basic} locally closed subset
$\Kbmo{X}{\sigma}$.  Thus to understand the
irreducible components of $\Kbm{X}{\tau}$ it suffices to understand
the irreducible components which pass through the general point of a
basic locus $\Kbmo{X}{\sigma}$.

\

\textbf{Convention}  Suppose that we have a contraction
$\phi:\sigma\rightarrow \tau$.  There is an induced morphism
\begin{equation}
\Kbmo{X}{\phi}:\Kbmo{X}{\sigma}\rightarrow \Kbm{X}{\tau}
\end{equation}
which is unramified with locally closed image.  We will speak of
$\Kbmo{X}{\sigma}$ as though it is a substack of $\Kbm{X}{\tau}$.
Thus given an irreducible component $M\subset \Kbmo{X}{\tau}$ and an
irreducible component $N\subset \Kbmo{X}{\sigma}$ we will say
$N\subset \overline{M}$ to mean that the image of $N$ is contained in
$\overline{M}$.  

\

The basic lemma is the following easy version of Mori's bend-and-break
lemma.  

\begin{lem}~\label{lem-bb}
Let $e>0$.
There is no complete curve contained in a fiber of the evaluation
morphism 
\begin{equation}
\text{ev}_{f_1,f_2}:\Kbmo{X}{\tau_2(e)}\rightarrow X\times X.
\end{equation}
\end{lem}

\begin{proof}
  Suppose that $C$ is a complete curve and $\zeta:C\rightarrow
  \Kbmo{X}{\tau_2(e)}$ has image in a fiber of $\text{ev}_{f_1,f_2}$,
  say the fiber over the pair $(p_1,p_2)\in X\times X$.  Denote the
  family $\zeta$ of strict $\tau_2(e)$-maps by $(\pi:\Sigma\rightarrow
  C,h,(q_1,q_2))$.  Let $H\subset X$ be a hyperplane section
  containing neither $p_1$ nor $p_2$.  Define $C'=\Sigma\times_X H$,
  so $C'$ is a finite ramified cover of $C$.  Let $B$ be the
  normalization of an irreducible component of $C'$ which dominates
  $C$.  The base-change $\Sigma \times_C B$ now admits the two
  sections $q_1,q_2$ as well as a third section $q_3$ which is
  everywhere disjoint from both $q_1$ and $q_2$.  Any $\PP^1$-bundle
  with three everywhere disjoint sections is isomorphic to
  $\PP^1\times B$ and the three sections are constant sections
  $\{0\}\times B, \{1\}\times B, \{\infty\}\times B$.  But now the
  morphism $h:\Sigma\times_C B\rightarrow X$ contracts the sections
  $q_1$ and $q_2$.  By the rigidity lemma, \cite[p. 43]{AVar}, we
  conclude that $h$ factors through the projection $\Sigma\times_C
  B\cong \PP^1\times B \rightarrow \PP^1$.  So we conclude that
  $C\rightarrow \Kbmo{X}{\tau}$ is a constant map.
\end{proof}

\begin{cor}~\label{cor-bb}
Let $M\subset \Kbm{X}{\tau_2(e)}$ be an irreducible, closed substack and
suppose that the fibers of $\text{ev}_{f_1,f_2}:M\rightarrow X\times
X$ have dimension at least $1$.  Then $M\cap
\lt(\Kbm{X}{\tau_2(e)}- \Kbmo{X}{\tau_2(e)}\rt)$ is either all of $M$
or contains an irreducible component with codimension $1$ in $M$.
\end{cor}

\begin{proof}
Suppose that $M$ intersects $\Kbmo{X}{\tau_2(e)}$.  Define
$I\subset X\times X$ to be the image $I=\text{ev}_{f_1,f_2}(M)$.  In
order to prove that 
\begin{equation}
\partial M :=
M\cap\lt(\Kbm{X}{\tau_2(e)}-\Kbmo{X}{\tau_2(e)}\rt)
\end{equation}
has an irreducible component of codimension $1$ in $M$, it suffices to
prove that for every $(p_1,p_2)\in I$, $\text{ev}^{-1}(p_1,p_2)\cap
\partial M\subset \text{ev}^{-1}(p_1,p_2)\cap M$ has codimension $1$.  
Suppose that it has codimension at least $2$.  Then the coarse moduli
space $|\text{ev}^{-1}(p_1,p_2)\cap \partial M|\subset
|\text{ev}^{-1}(p_1,p_2) \cap M|$ has codimension at least $2$.  Since
the coarse moduli spaces are proper varieties, we can find a complete
curve $C$ in $|\text{ev}^{-1}(p_1,p_2)\cap M|$ which does not
intersect $|\text{ev}^{-1}(p_1,p_2)\cap \partial M|$.  Since
$\Kbm{X}{\tau_2(e)}$ is a Deligne-Mumford stack, there exists a finite
ramified cover $C'\rightarrow C$ such that $C'\rightarrow
|\Kbm{X}{\tau_2(e)}|$ factors through $\Kbm{X}{\tau_2(e)}\rightarrow
|\Kbm{X}{\tau_2(e)}|$.  But then by lemma~\ref{lem-bb}, we conclude
that $C'\rightarrow |\Kbm{X}{\tau_2(e)}|$ is constant, which
contradicts the construction of $C$.  Therefore $\partial M\subset M$
has codimension $1$.
\end{proof}

The main application is the following:

\begin{prop}~\label{prop-bbflat}
Suppose that $X\subset \PP^N$ is a complete intersection.  
Suppose that $\fe(X,\tau_1(e),f_1)$ holds for every $e <
E$
and suppose that every irreducible component of $\Kbmo{X}{\tau_1(E)}$
has dimension at least $2\text{dim}(X)$.  Then
$\fe(X,\tau_1(E),f_1)$ holds as well.  
\end{prop}

\begin{proof}
By lemma~\ref{lem-propeq}, to prove that
$\fe(X,\tau_1(E),f_1)$ holds, it suffices to prove that
$\cE(X,\tau_1(E),f_1)$ holds.  Let $\zeta\in
\Kbmo{X}{\tau_1(e)}$ be a point, denote $p=\text{ev}_{f_1}(\zeta)$,
and let $M\subset \text{ev}_{f_1}^{-1}(p)$ be an irreducible component
which contains $p$.  We need to prove that
$\text{dim}(M)=\text{dim}(X,\tau_1(e)) -\text{dim}(X)$.  

\

Now consider the forgetful morphism
\begin{equation}
\Kbmo{X}{a}:\Kbmo{X}{\tau_2(e)}\rightarrow \Kbmo{X}{\tau_1(e)}.
\end{equation}  
This is a smooth surjective morphism: let $N\subset
\Kbmo{X}{\tau_2(e)}$ be the preimage of $M$.  Then
$\text{dim}(N)=\text{dim}(M)+1$, thus we have to prove that 
\begin{equation}
\text{dim}(N) = \text{dim}(X,\tau_1(e))+1 - \text{dim}(X).
\end{equation}
It suffices to prove that the general fiber of $\text{ev}_{f_2}:N\rightarrow
X$ has dimension at most $\text{dim}(X,\tau_1(e))+1 - 2\text{dim}(X)$
(since we already know the dimension is at least this large).

\

Choose any point $q\neq p$ in $\text{ev}_{f_2}(N)$ and consider $N_q
:= \text{ev}_{f_2}^{-1}(q)\cap N$.  
By assumption, $\text{dim}(X,\tau_1(e))+1-2\text{dim}(X) \geq 1$, so
$\text{dim}(N_q) \geq 1$.  
Define $\overline{N}\subset \Kbm{X}{\tau_2(e)}$ to be the closure of
$N$.
By corollary~\ref{cor-bb} we conclude that $\partial
\overline{N}\subset \overline{N}$ has codimension $1$.  In other
words, there is a stable $A$-graph $\sigma'\not\cong\tau_2(e)$ whose
canonical 
contraction is $\phi':\sigma'\rightarrow \tau_2(e)$ and such that
$\overline{N}\cap \text{Image}(\Kbmo{X}{\sigma'})\subset \overline{N}$
has an irreducible component of codimension $1$.  

\

Now there is precisely one stable $A$-graph $\sigma_0\not\cong
\tau_2(e)$ whose stabilization after removing $f_2$ equals
$\tau_1(e)$, namely:

\begin{picture}(250,150)
\put(50,50){\makebox(200,100){\includegraphics{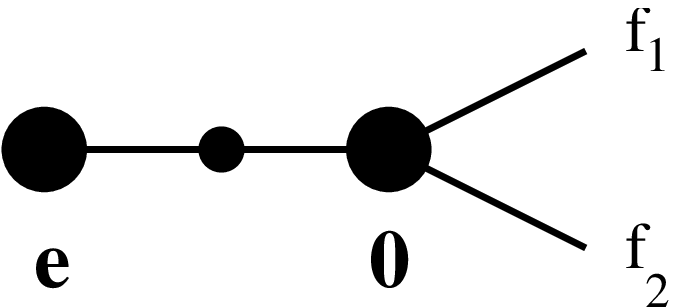}}}
\put(150,20){\makebox(0,0)[b]{Diagram 2}}
\end{picture}


By the assumption that $q\neq p$, the point $(q,p)\not\in
\text{ev}_{f_1,f_2}(\Kbmo{X}{\sigma_0})$.  So $\sigma'$ is
not $\sigma_0$.  We conclude that the image of
$\Kbmo{X}{\sigma'}$ in $\Kbm{X}{\tau_1(e)}$ (under the map which
stably removes $f_2$) is again a boundary component $\Kbm{X}{\sigma}$
for some $\phi:\sigma\rightarrow \tau_1(e)$ (not the identity).  Therefore
$\overline{M}\cap\text{Image}(\Kbmo{X}{\sigma})\subset\overline{M}$
is a locally closed substack such that some irreducible component has
codimension one in $\overline{M}$.  

\

Since $\sigma\not\cong\tau_1(E)$, we have $E(\sigma)< E$.  By our
assumption and by proposition~\ref{prop-red}, we conclude that
$\fe(X,\sigma,f_1)$ holds.  In particular,
$\overline{M}\cap\text{Image}(\Kbmo{X}{\sigma})$ has dimension at most
$\text{dim}(X,\sigma)-\text{dim}(X)$.  So the
dimension of $\overline{M}$ is at most
$\text{dim}(X,\sigma)+1-\text{dim}(X)$.  Since
$\text{dim}(X,\sigma)+1\leq \text{dim}(X,\tau_1(e))$, we conclude
that $\text{dim}(M)\leq \text{dim}(X,\tau_1(e))-\text{dim}(X)$.  So
$\fe(X,\tau_1(e),f_1)$ holds.  
\end{proof}

\textbf{Remark}  The condition that $X\subset \PP^N$ be a complete
intersection is essentially superfluous.  If instead we had worked
throughout with the property $\cE(X,\tau,f)$ rather than
$\fe(X,\tau,f)$ (which is a little trickier), then the
argument above proves the analogous result without this condition on
$X$.  

\begin{defn}~\label{defn-thrdeg}
If $X\subset \PP^N$ is a complete intersection of hypersurfaces of
degrees $d_1,d_2,\dots,d_r$, define the \emph{threshold degree} to be
\begin{equation}
E(X)=E(N,(d_1,\dots,d_r)) = \lt\lfloor \frac{N+2-r}{N+1-(d_1+\dots+d_r)}
\rt\rfloor.
\end{equation}
In particular, if $X\subset \PP^N$ is a hypersurface of degree $d <
\frac{N+1}{2}$ then $E(X) = 1$.
\end{defn}

\begin{cor}~\label{cor-bbflat}
Suppose that $X\subset \PP^N$ is a complete intersection of
hypersurfaces of degrees $d_1,d_2,\dots,d_r$.  If
$\fe(X,\tau_1(e),f_1)$ holds for each $1\leq e \leq E(X)$,
then for every stable $A$-graph $\tau$ and flag
$f\in\text{Flag}(\tau)$, $\fe(X,\tau,f)$ holds.
In particular, if $X\subset \PP^N$ is a hypersurface of degree
$d<\frac{N+1}{2}$ then it suffices to prove
$\fe(X,\tau_1(1),f_1)$.  
\end{cor}

\begin{proof}
By proposition~\ref{prop-red}, to prove that
$\fe(X,\tau,f)$ always holds it suffices to prove that
$\fe(X,\tau_1(e),f_1)$ always holds (for $e>0$).  Now for
$e>E(X)$ we have by lemma~\ref{lem-lbdim} that every irreducible
component of $\Kbmo{X}{\tau_1(e)}$ has dimension at least
\begin{equation}
\begin{array}{c}
 \lt(N+1-(d_1+\dots +d_r)\rt)e + (N-r-3) + 1 \geq \\ 
\lt(N+1-(d_1+\dots+ d_r)\rt)\frac{N+2-r}{N+1-(d_1+\dots+d_r)} +
(N-r-3) + 1,
\end{array}
\end{equation}
which, of course, is just $2\text{dim}(X)$.  
So by proposition~\ref{prop-bbflat} and induction, to
prove $\fe(X,\tau_1(e),f_1)$ for all $e$, it suffices to
prove $\fe(X,\tau_1(e),f_1)$ in the cases $e=1,\dots,E(X)$.
\end{proof}

\begin{cor}~\label{cor-bbflat2}
With the same hypotheses as in corollary~\ref{cor-bbflat}, suppose
that $\fe(X,\tau_1(e),f_1)$ holds for each $1\leq e\leq
E(X)$.  Then for every stable $A$-graph $\tau$, $\Kbmo{X}{\tau}$ has
pure dimension $\text{dim}(X,\tau)$.  
\end{cor}

\begin{proof}
If $E(\tau)=0$, then $\Kbmo{X}{\tau}=\Kbmo{*}{\tau}\times X$ and the
result follows from ~\cite[proposition 7.4]{BM}.  
Suppose that $E(\tau)>0$ and let
$v\in\text{Vertex}(\tau)$ be such that $\beta(v)>0$.  Define
$a:\tau'\hookleftarrow \tau$ to be the combinatorial morphism which
adds a new tail $f$ to $v$.  Then $\Kbmo{X}{a}$ is smooth and
surjective of relative
dimension $1$.  By corollary~\ref{cor-bbflat},
$\fe(X,\tau',f)$ holds.  In particular, $\Kbmo{X}{\tau'}$
has pure dimension $\text{dim}(X,\tau')$.  It follows that
$\Kbmo{X}{\tau}$ has pure dimension $\text{dim}(X,\tau')-1 =
\text{dim}(X,\tau)$.  
\end{proof}

A second application of the proof of proposition~\ref{prop-bbflat} is
the following:

\begin{prop}~\label{prop-bbspec}
Suppose that $X\subset \PP^N$ is a complete intersection.  Suppose
that $\fe(X,\tau_1(e),f_1)$ holds for every $e< E$ and
suppose that every irreducible component of $\Kbmo{X}{\tau_1(E)}$ has
dimension at least $2\text{dim}(X)$.  Then for every irreducible
component $M\subset\Kbmo{X}{\tau_0(E)}$ there is a graph
$\sigma=\tau_{0,0}(i,j)$, $0<i,j$ and $i+j=E$, and an irreducible
component $N\subset \Kbmo{X}{\tau_{0,0}(i,j)}$ such that $N\subset
\overline{M}$ is a codimension $1$ subvariety.
\end{prop}

\begin{proof} 
Let $M'\subset \Kbmo{X}{\tau_1(e)}$ be the irreducible component which
dominates $M$.  By the proof of proposition~\ref{prop-bbflat} 
there is a graph $\sigma\not\cong\tau_1(e)$ with canonical
contraction $\phi:\sigma\rightarrow \tau_1(e)$ such that
$\overline{M'}\cap \Kbmo{X}{\sigma}\subset \overline{M}$ has
codimension $1$.  Now $\Kbmo{X}{\sigma}$ has dimension
$\text{dim}(X,\sigma)$ and by lemma~\ref{lem-lbdim} $\overline{M'}$ has
dimension at least 
$\text{dim}(X,\tau_1(e))$.  But $\text{dim}(X,\sigma) =
\text{dim}(X,\tau_1(e)) - \#\text{Edge}(\sigma)$.  So 
$\sigma$ has exactly one edge, i.e. $\sigma=\tau_{1,0}(i,j)$ for some
$i,j$ with $i+j=E$.  By stability $i,j>0$.  Moreover,
$\overline{M'}\cap \Kbmo{X}{\sigma}$ has dimension
$\text{dim}(X,\sigma)$ so there is an irreducible
component $N'\subset \Kbmo{X}{\sigma}$ such that $N'\subset
\overline{M'}$.  

\

Since $\Kbmo{X}{\tau_{1,0}(i,j)}\rightarrow
\Kbmo{X}{\tau_{0,0}(i,j)}$ is smooth of relative dimension $1$, there
is an irreducible component $N\subset \Kbmo{X}{\tau_{0,0}(i,j)}$ such
that $N'$ is the preimage of $N$.  Thus $N\subset \overline{M}$ is a
codimension $1$ subvariety.
\end{proof}

\begin{defn}~\label{defn-basic}
Let $X\subset \PP^N$ be a complete intersection with
threshold degree
$E(X)=E$.  A stable $A$-graph $\tau$ is \emph{basic} for $X$ if its
maximal component degree $E(\tau)$ satisfies
$E(\tau)\leq E(X)$.  
\end{defn}

\begin{defn}~\label{defn-nice}
For an $A$-graph $\tau$, define its \emph{degree 0
subgraph} to be the maximal subgraph $\tau\hookleftarrow \tau^0$ such
that $E(\tau^0)=0$, i.e. $\tau^0$ is the (possibly disconnected)
subgraph of $\tau$ with 
\begin{equation}
\text{Vertex}(\tau^0)=\{v\in\text{Vertex}(\tau)| \beta(v)=0\}
\end{equation}
and 
\begin{equation}
\text{Flag}(\tau^0)=\{f\in\text{Flag}(\tau)| \partial
f\in\text{Vertex}(\tau^0)\}.
\end{equation}
A contraction of $A$-graphs $\phi:\sigma\rightarrow \tau$ is
\emph{nice} if $\phi$ induces an isomorphism $\sigma^0\cong\tau^0$.
\end{defn}

\begin{thm}~\label{thm-spec}
Let $X\subset \PP^N$ be a complete intersection, $\tau$ a stable
$A$-graph and $M\subset \Kbmo{X}{\tau}$ an irreducible component.
Suppose that $\fe(X,\tau_1(e),f_1)$ holds for each $1\leq
e\leq E(X)$.  Then there exists a nice contraction
$\phi:\sigma\rightarrow \tau$ and an irreducible component $N\subset
\Kbmo{X}{\sigma}$ such that $\sigma$ is basic and such that $N\subset
\overline{M}$.
\end{thm} 

\begin{proof}
We will prove this by induction on the maximal component degree
$E(\tau)$.  If $E(\tau) \leq E(X)$, then we can take $\phi$ to be the
identity $\tau\rightarrow \tau$ and $N=M$.  

\

Suppose that $E > E(X)$.  By way of induction, 
assume the theorem is proved for all
graphs $\tau$ with $E(\tau) < E$.  We will deduce the theorem for
graphs with $E(\tau)=E$ by induction on $\#\text{Vertex}(\tau)$, and
thus establish the theorem by induction on $E(\tau)$.

\

First we consider the case that $\#\text{Vertex}(\tau)=1$,
i.e. $\tau=\tau_r(E)$ 
for some $r$.  Define $a:\tau_r(E)\hookleftarrow \tau_0(E)$ be the
combinatorial morphism which strips the tails from $\tau_r(E)$.  
Then $\Kbmo{X}{a}:\Kbmo{X}{\tau_r(E)}\rightarrow
\Kbmo{X}{\tau_0(E)}$ is smooth, surjective with connected fibers of
dimension $r$.  So we 
conclude that $M$ is the preimage of an irreducible component
$M'\subset \Kbmo{X}{\tau_0(E)}$.  Now by proposition~\ref{prop-bbspec}
there is a nice contraction $\psi:\rho\rightarrow \tau_0(E)$ and an
irreducible component $L\subset\Kbmo{X}{\rho}$ such that $L\subset
\overline{M}'$.  As $E(\rho)<E$, by assumption there exists
a nice contraction $\phi':\sigma'\rightarrow \rho$ such that $\sigma'$ is
basic and there exists an irreducible component $N'\subset
\Kbmo{X}{\sigma}$ such that $N'\subset \overline{L}$.  Therefore
$\psi\circ\phi':\sigma'\rightarrow \tau_0(E)$ is a nice contraction and
$N'\subset \overline{M'}$.  Now let $v\in\text{Vertex}(\sigma')$ be
any vertex and let $b:\sigma\hookleftarrow \sigma'$ be the graph
obtained by attaching $r$ tails to $v$.  Note that $\Kbmo{X}{b}$ is
smooth, surjective with connected fibers of dimension $r$.  
Let $\phi:\sigma\rightarrow
\tau_r(E)$ be the contraction obtained from $\psi\circ\phi'$ by
sending the $r$ tails of $\sigma$ to the $r$ tails of $\tau_r(E)$.
Then $\phi$ is a nice contraction and $\sigma$ is basic.  Define
$N\subset \Kbmo{X}{\sigma}$ to be the preimage of $N'$ under
$\Kbmo{X}{a}$.  
The morphism $\Kbm{X}{a}$ is compatible with $\Kbmo{X}{b}$,
i.e. $\Kbm{X}{a}\circ \Kbmo{X}{\phi} =
\Kbmo{X}{\psi\circ\phi'}\circ\Kbmo{X}{b}$, and $\Kbm{X}{a}$ is smooth
along the image of $\Kbmo{X}{\phi}$.  Thus we conclude that $N\subset
\overline{M}$, the theorem is proved for $M$.

\

Now we consider the general case.  
For each graph $\tau$, define
\begin{equation}
l(\tau):= \#\text{Vertex}(\tau).
\end{equation}
When $l=1$, we have the case in the last paragraph.  Suppose $l>1$
and, by way of induction on $l$, 
suppose that for all stable $A$-graphs $\tau$ with
$E(\tau)=E$ and with $l(\tau) < l$, the theorem is proved.
Suppose that $\tau$ is a stable $A$-graph with $E(\tau)=E$
and $l(\tau)=l$.  Let $\{f_1,f_2\}$ be an edge of
$\tau$.  Define $a_1:\tau\hookleftarrow \tau_1$ and
$a_2:\tau\hookleftarrow \tau_2$ be the two subgraphs obtained by
breaking the edge (see diagram 1).

\

Now the
morphism $\Kbmo{X}{a_1}:\Kbmo{X}{\tau}\rightarrow \Kbmo{X}{\tau_1}$ is
the composition of an open immersion and the projection of the fiber
product
\begin{equation}
\Kbmo{X}{\tau_1}\times_{\text{ev}_{f_1},X,\text{ev}_{f_2}}\Kbmo{X}{\tau_2}\rightarrow \Kbmo{X}{\tau_1}.
\end{equation}
By proposition~\ref{prop-bbflat} the morphism $\text{ev}_{f_2}$ is
flat.  Therefore $\Kbmo{X}{a_1}$ is flat.  So $M$ dominates an
irreducible component $M_1$ of $\Kbmo{X}{a_1}$.  Since $l(\tau_1) =
l(\tau) - l(\tau_2)$ and $l(\tau_2) > 0$ by the assumption that
$E(\tau_2) = E$, we have $l(\tau_1)<l(\tau)$. So by the induction
assumption, there
exists a nice contraction $\phi_1:\rho_1\rightarrow \tau_1$ and an
irreducible component $L_1\subset \Kbmo{X}{\rho_1}$ such that
$\rho_1$ is basic and such that $L_1\subset \overline{M_1}$.  

\

Since $\Kbm{X}{a_1}$ is proper, there
exists an irreducible subvariety $L\subset \overline{M}$ such that
$\Kbm{X}{a_1}(L)=L_1$ and such that the fiber dimension of
$L\rightarrow L_1$ is at least the fiber dimension of
$\Kbmo{X}{a_1}$.  Up to replacing $L$ by an open subset, we may
suppose $L$ is contained in one of the locally
closed substacks $\Kbmo{X}{\rho}$.  Since $\Kbm{X}{a_1}$ maps $L$ into
$\Kbmo{X}{\rho_1}$, we must have that $\rho$ is glued from
$\rho_1$ and a graph $\rho_2$ by making an edge out of $f_1$ and
$f_2$.  Moreover $\rho_2$ must contract to $\tau_2$.  But then the
dimension of $L$ is at most 
\begin{equation}
\begin{array}{c}
\text{dim}(X,\rho_1) + \text{dim}(X,\rho_2) -
\text{dim}(X) = \\
\text{dim}(X,\rho_1) +
\lt(\text{dim}(X,\tau)-\text{dim}(X,\tau_1) \rt)
-\lt(\#\text{Edge}(\rho_2) - \#\text{Edge}(\tau_2)\rt).
\end{array}
\end{equation}
Thus we conclude that $\#\text{Edge}(\rho_2)=\#\text{Edge}(\tau_2)$,
i.e. $\rho_2=\tau_2$.  So $\psi:\rho\rightarrow \tau$ is a nice
contraction and $L\subset \Kbmo{X}{\rho}$ is an irreducible component
such that $L\subset \overline{M}$.  Moreover, $l(\rho) = l(\rho_2) =
l(\tau_2) < l(\tau)$.  By the induction assumption, there exists a
nice contraction $\phi:\sigma\rightarrow \rho$ and an irreducible
component $N\subset \Kbmo{X}{\sigma}$ such that $\sigma$ is basic and
such that $N\subset \overline{L}$.  But then
$\psi\circ\phi:\sigma\rightarrow \tau$ is nice and
$N\subset\overline{M}$, i.e. the theorem is proved for $M$.  So the
theorem is proved by induction on $E$ and $l$.  
\end{proof}  

The previous theorem suggests a strategy for proving that any given
$\Kbm{X}{\tau}$ is irreducible:  
\begin{enumerate}
\item Determine all nice contractions $\phi:\sigma\rightarrow \tau$
such that $\sigma$ is basic.
\item Determine all irreducible components $N\subset\Kbmo{X}{\sigma}$.
\item Show that for each $N$, there is a unique irreducible component
$M(N)\subset \Kbm{X}{\tau}$ which contains $N$.
\item Prove that all of the putative irreducible components $M(N)$ are
actually equal.
\end{enumerate}
The first step (1) is a combinatorial problem.  
The simplest case for (2) is when $\Kbmo{X}{\sigma}$ is itself
irreducible for each basic $\sigma$.
One can try to prove (3) by a deformation theory argument; if one
proves that the general point of $M$ is a smooth point of the stack
$\Kbm{X}{\tau}$, then it follows that there is a unique irreducible
component $M(N)$ which contains $N$.
The proof of (4) reduces to a combinatorial argument: there is a
combinatorially-defined equivalence relation on nice contractions such
that equivalent contractions give rise to the same irreducible
component $M(N)$.  Thus we are reduced to proving that all nice
contractions $\sigma \rightarrow \tau$ are equivalent for this
equivalence relation. 

\section{Properties of evaluation morphisms}~\label{sec-basic}

In the last two sections we investigated when an evaluation morphism
$\text{ev}_f:\Kbmo{X}{\tau}\rightarrow X$ is flat of the expected
dimension.  In this section we also investigate when the general fiber
is irreducible, and when the morphism $\text{ev}_f$ is
\emph{unobstructed} at a general point of $\Kbmo{X}{\tau}$.  By the
same techniques in the last two sections, we reduce these properties
for a general basic $A$-graph $\tau$ to the property for $A$-graphs of
the form $\tau_1(e)$ with $e=1,\dots,E(X)$.  Using this result, we
carry out steps $(2)$ and $(3)$ of the strategy of proof in the
previous sections.

\

The new property of evaluation morphisms we want to consider is the
following. 

\begin{defn}  Suppose $X\subset \PP^N$ is a smooth subvariety, $\tau$
  is a stable $A$-graph and $f\in\text{Flag}(\tau)$.   We say that
  $\cB(X,\tau,f)$ holds if we have
\begin{enumerate}
\item $\fe(X,\tau,f)$ holds,
\item the general fiber of $\text{ev}_{f}$ is geometrically
  irreducible,
\item in $\Kbmo{X}{\tau}$ there is a point $\lt[h:C\rightarrow
  X\rt]$ which is \emph{free}, i.e. $h^*T_X$ is generated by global
  sections.
\end{enumerate}
\end{defn}

The difficult item to check is still $(1)$.  We return to this point
at the end of this section.  First we give a reformulation of
item $(3)$ of the definition above.

\begin{lem}~\label{lem-Breform}
Suppose that $(C,q_{f_i})$ is a genus $0$ prestable curve with dual
graph $\tau$, and suppose that $E$ is a locally free sheaf on $C$.
The following are equivalent:
\begin{enumerate}
\item $E$ is generated by global sections.
\item For each irreducible component $C_v$ of $C$, the restriction
  $E_v$ of $E$ to $C_v$ is generated by global sections.
\item For each irreducible component $C_v$ of $C$ and each point $q\in
  C_v$, we have $H^1(C_v,E_v(-q))= 0$.
\end{enumerate}
\end{lem}

\begin{proof}
Clearly $(1)$ implies $(2)$.  By Grothendieck's lemma~\cite[exercise
V.2.4]{H}, $E_v$ splits as a direct sum of line bundles
$L_1\oplus\dots \oplus L_r$.  If $(2)$ is satisfied, then each $L_v$
is generated by global sections, i.e. if we identify $C_v$ with
$\PP^1$, then $L_i=\OO_{\PP^1}(a_i)$ with $a_i\geq 0$.  Since
$H^1(\PP^1,\OO_{\PP^1}(a-1)) =0$ for $a\geq 0$, we conclude that
$H^1(C_v,E_v(-q)) = 0$.  So $(2)$ implies $(3)$.  

\

Finally suppose that $(3)$ is satisfied.  We shall prove that $E$ is
generated by global sections by induction on the number of vertices of
$\tau$, i.e. on the number of irreducible components of $C$.  If $C$
has a single irreducible component, then 
$C$ is isomorphic to $\PP^1$.  By Grothendieck's lemma we know
$E=\OO_{\PP^1}(a_1)\oplus \dots \oplus \OO_{\PP^1}(a_r)$ for some
integers $a_i$.  Since $H^1(\PP^1,E(-q))=0$, we conclude that each
$a-1\geq -1$, i.e $a\geq 0$.  So $E$ is generated by global sections.

\

Now suppose that $\tau$ has more than one vertex and let $v_1$ be any leaf
of $\tau$, i.e. $v_1$ is adjacent to exactly one other vertex.
Let $i_1:C_1\rightarrow C$ be the irreducible component associated to
$v_1$, let $i_2:C_2\rightarrow C$ be the union of all the other
irreducible components of $C$ and let $q$ be the unique point of
intersection of $C_1$ and $C_2$.  Let $E_1$ denote the restriction of
$E$ to $C_1$ and let $E_2$ denote the restriction of $E$ to $C_2$.  By
the induction assumption, we may assume that $E_2$ is generated by
global sections.  But now we have an exact sequence of sheaves on $C$:
\begin{equation}
0\rightarrow (i_1)_*(E_1(-q)) \rightarrow E \rightarrow (i_2)_*(E_2)
\rightarrow 0.
\end{equation}
The obstruction to lifting the global sections of $E_2$ to global
sections of $E$ is an element of $H^1(C_1,E_1(-q))$, which is zero by
assumption.  So every global section of $E_2$ is the restriction of a
global section of $E$.  Thus the locus where $E$ isn't generated by
global sections (i.e. the cokernel of the morphism $H^0(C,E)\otimes_\CC
\OO_E\rightarrow E$) is a closed subset of $C_1 - C_2$.  Since $\tau$ has more
than one vertex, we can find a second leaf $v_2$ of $\tau$.  Repeating
the argument with $v_2$ we conclude that $E$ is generated by global
sections.
\end{proof}

\begin{lem}~\label{lem-Breform2}
With the same notation as in lemma~\ref{lem-Breform}, if $E$ satisfies
any of the three equivalent conditions above, and if $p\in C$ is any
smooth point, then $H^1(C,E(-p))=0$.
\end{lem}

\begin{proof}
If $C$ has a single irreducible component, this follows from the
equivalent condition $(3)$ in lemma~\ref{lem-Breform}.  
Suppose that
$C$ has $l>1$ irreducible components.  By way of induction, suppose
that the lemma has been proved for all curves with fewer than $l$
irreducible components.
We can find a leaf
$v_1$ of $C$ such that $p$ is not contained in the corresponding
irreducible component $C_2$.  
Let $C_2$, $i_1$, $i_2$, and $q$ be as in the proof of
lemma~\ref{lem-Breform}.  Then we have a short exact sequence:
\begin{equation}
0 \rightarrow (i_1)_*(E_1(-q)) \rightarrow E(-p) \rightarrow
(i_2)_*(E_2(-p)) \rightarrow 0.
\end{equation}
By the induction assumption, both $H^1(C_1,E_1(-q))=0$ and
$H^1(C_2,E_2(-p))=0$.  So by the long exact sequence in cohomology
associated to the short exact sequence above, we conclude that
$H^1(C,E(-p))=0$.  So the lemma is proved by induction on $l$.
\end{proof}

\begin{lem}~\label{lem-Blines}
Suppose $n>2$ and $X\subset \PP^n$ is a general hypersurface of degree
$d<\frac{n+1}{2}$.  Then $\cB(X,\tau_1(1),f_1)$ holds.
\end{lem}

\begin{proof}
By theorem~\ref{thm-thm1}, for general $X$ we have that
$\fe(X,\tau_1(1),f_1)$ holds, i.e. we have (1).  
By lemma~\ref{lem-dfmlines}, for general
$X$ there is a free line on $X$, i.e. we have (3).  

\

By theorem~\ref{lem-smlines}, for
general $X$, $F_{0,1}(X)$ is smooth.  Thus by generic smoothness, the
general fiber of $\text{ev}_{f_1}$ is smooth.  By lemma~\ref{lem-lem1},
the general fiber of $\text{ev}_{f_1}$ is a complete intersection in
$\PP^{n-1}$ of dimension $n-d-1 > 1$.  Thus the general fiber is
geometrically 
connected by repeated application of ~\cite[corollaire 3.5,
exp. XII]{SGA2}.  Since a smooth, geometrically connected scheme is
geometrically irreducible, we have (2).
\end{proof}

The main result of this section is the following:

\begin{prop}~\label{prop-basic}
Suppose $X\subset \PP^N$ is a smooth subvariety which satisfies
$\cB(X,\tau_1(e),f_1)$ for $e=1,\dots, E$.  
Let $\tau$ be an $A$-graph such that $E(\tau)\leq E$.  
Then we have the following:
\begin{enumerate}
\item For each $f\in\text{Flag}(\tau)$, we have $\cB(X,\tau,f)$.
\item $\Kbmo{X}{\tau}$ is an irreducible stack.
\end{enumerate}
\end{prop}

\begin{proof}
Both statements are trivial in case $\tau$ is empty, so assume $\tau$
is nonempty.  Observe that $(1)$ implies $(2)$: given
$v\in\text{Vertex}(\tau)$, define a new $A$-graph $\tau'$ and a
combinatorial morphism $\alpha:\tau'\hookleftarrow \tau$ which
attaches 
a new
flag $f'$ to $\tau$ at $v$.  Then
$\Kbmo{X}{\alpha}:\Kbmo{X}{\tau'}\rightarrow \Kbmo{X}{\tau}$ is
smooth, surjective with geometrically irreducible fibers.  So
$\Kbmo{X}{\tau'}$ is irreducible iff $\Kbmo{X}{\tau}$ is irreducible.
By $(1)$, $\text{ev}_{f'}:\Kbmo{X}{\tau'}\rightarrow X$ is flat and the
general fiber is geometrically irreducible.  Since $X$ is irreducible,
it follows that $\Kbmo{X}{\tau'}$ is irreducible.  So it remains to
prove $(1)$.

\

First of all, suppose that $\beta(\tau)=0$.  Let
$\alpha:\tau\hookleftarrow \emptyset$ be the unique morphism.  Then
$\text{ev}_f$ coincides with
$\Kbmo{X}{\alpha}:\Kbmo{X}{\tau}\rightarrow \Kbmo{X}{\emptyset}=X$.
Thus $(1)$ follows from lemma~\ref{lem-isog}.  So we are reduced to the
case that $\beta(\tau)> 0$.

\

We prove $(1)$ by induction on the number of vertices of $\tau$.
Suppose that $\tau$ has a single vertex, i.e. $\tau = \tau_r(e)$ for
some $r>0$.  Let $\alpha:\tau_r(e)\hookleftarrow \tau_1(e)$ be the
unique combinatorial morphism which maps
$f_1\in\text{Flag}(\tau_1(e))$ to $f\in\text{Flag}(\tau_r(e))$.  Then
$\text{ev}_f$ factors as the composition
\begin{equation}
\begin{CD}
\Kbmo{X}{\tau_r(e)} @>\Kbmo{X}{\alpha} >> \Kbmo{X}{\tau_1(e)} @>
\text{ev}_{f_1} >> X.
\end{CD}
\end{equation}
Of course $\Kbmo{X}{\alpha}$ is an open immersion into the
$(r-1)$-fold fiber product of the universal curve, so
$\Kbmo{X}{\alpha}$ is smooth with geometrically irreducible fibers.
And by $\cB$, $\text{ev}_{f_1}$ is flat and the general fiber is
geometrically irreducible.  Therefore the composition is flat and the
general fiber is geometrically irreducible, i.e. items $(1)$ and $(2)$
of condition $\cB(X,\tau_r(e),f)$ hold.  The condition that
$h^*T_X$ is generated by global sections is independent of the number
of vertices, so item $(3)$ of $\cB(X,\tau_r(e),f)$ holds as well,
i.e. $\cB(X,\tau_r(e),f)$ holds.

\

Now suppose that there is more than one vertex, say
$\#\text{Vertex}(\tau) = l > 1$.    
By way of assumption, suppose that $\cB(X,\sigma,f)$ has been proved for all 
$A$-graphs $\sigma$ such that $\#\text{Vertex}(\sigma) < l$.  
Let $\{f_1,f_2\}$ be
any edge and consider the subgraphs $\alpha_1:\tau\hookleftarrow
\tau_1$ and $\alpha_2:\tau\hookleftarrow \tau_2$ as in diagram $1$.
Without loss of generality, suppose that $f$ is in $\tau_1$.  Then 
$\text{ev}_f$ factors as the composition:
\begin{equation}
\begin{CD}
\Kbmo{X}{\tau} @> \Kbmo{X}{\alpha_1} >> \Kbmo{X}{\tau_1} @>
\text{ev}_f >> X.
\end{CD}
\end{equation}
Now $\Kbmo{X}{\alpha_1}$ factors as the composition of an open
immersion and the projection
\begin{equation}
\pi_1: \Kbmo{X}{\tau_1}\times_{\text{ev}_{f_1},X,\text{ev}_{f_2}}
\Kbmo{X}{\tau_2} \rightarrow \Kbmo{X}{\tau_1}.
\end{equation}
Since $\#\text{Vertex}(\tau_i) < l$ for $i=1,2$, 
the induction assumption says
that $(2)$ holds for $\text{ev}_{f_i}:\Kbmo{X}{\tau_i} \rightarrow
X$.  Since $\text{ev}_{f_1}$ is open, the general fiber of $\pi_1$
dominates the general fiber of $\text{ev}_{f_2}$.  So $\pi_1$ is flat
and the general fiber is geometrically irreducible.  Thus the same is
true of $\Kbmo{X}{\alpha_1}$.  Since $\#\text{Vertex}(\tau_1) < l$, 
$\text{ev}_{f}:\Kbmo{X}{\tau_1}\rightarrow X$ is flat and the general
fiber is geometrically irreducible.  Thus the composition is flat and
the general fiber is geometrically irreducible, i.e. items $(1)$ and
$(2)$ of $\cB(X,\tau,f)$ hold.

\

Finally we consider item $(3)$ of $\cB(X,\tau,f)$.  Each of the two
projections $\Kbmo{X}{\alpha_1}$ and $\Kbmo{X}{\alpha_2}$ are
dominant.  By the induction assumption, for $i=1,2$ the set of points
in $\Kbmo{X}{\tau_1}$ which parametrize stable maps with $h^*T_X$
generated by global sections is an open, dense set $U_i$.  The preimage of
each $U_i$ in $\Kbmo{X}{\tau}$ is an open, dense set, and the
intersection of these two open, dense sets is an open, dense set.  For a
point in this intersection -- using the equivalent condition $(2)$ of
lemma~\ref{lem-Breform} -- we have that the restriction of $h^*T_X$ to
each irreducible component with vertex $v\in \tau_1$ is generated by
global sections, and also the restriction of $h^*T_X$ to each
irreducible component with vertex $v\in \tau_2$ is generated by global
sections.  So by the equivalent condition $(2)$ of
lemma~\ref{lem-Breform}, we conclude that $h^*T_X$ is generated by
global sections.  Thus item $(3)$ of $\cB(X,\tau,f)$ is satisfied.
This completes the proof that $\cB(X,\tau,f)$ holds, and the
proposition is proved by induction.
\end{proof}

Propsition~\ref{prop-basic} simplifies step 2 in the strategy of the last
section to checking that $\cB(X,\tau_1(e),f_1)$ holds for all
$e=1,\dots,E(X)$.  
Next we reduce step 3 to checking that $\cB(X,\tau_1(e),f_1)$ holds
for all $e=1,\dots,E(X)$.

\begin{prop}~\label{prop-Bdefm}
Suppose that $\tau$ is a stable $A$-graph, $f\in\text{Tail}(\tau)$ and
suppose that $\cB(X,\tau,f)$ holds.
Suppose that
$\alpha:\tau\rightarrow \sigma$ is a contraction. 
The morphism $\Kbmo{X}{\alpha}$ maps a
general point of $\Kbmo{X}{\tau}$ to a point in the smooth locus of
the morphism of $\text{ev}_f:\Kbm{X}{\sigma}\rightarrow X$.
\end{prop}

\begin{proof}
As in the proof of proposition~\ref{prop-basic}, the case that
$\beta(\sigma)=0$ follows from lemma~\ref{lem-isog}.  So we are reduced to
the case $\beta(\sigma) > 0$.  

\

Let $\mathfrak{M}(\sigma)$ denote the (non-separated) Artin stack of
prestable $\sigma$-curves as in the proof of lemma~\ref{lem-lbdim}.
There is a 
$1$-morphism $\Kbm{X}{\sigma}\rightarrow \mathfrak{M}(\sigma)$ given
by forgetting the map to $X$.  ``Remembering'' the map to $X$ gives an
open immersion of $\Kbm{X}{\sigma}$ into the relative scheme of morphisms
$\text{Mor}_{{\mathfrak M}(\sigma)}(\mathfrak{C}(\sigma),X \times
\mathfrak{M}(\sigma))$ (notation as in the proof of lemma~\ref{lem-lbdim}).
Since ${\mathfrak M}(\sigma)$ is smooth by ~\cite[prop. 2]{B}, to
prove that $\text{ev}_{f}$ is smooth at a point, it suffices to prove
the following morphism is smooth at this point:
\begin{equation}
(\pi,\text{ev}_f):\Kbm{X}{\sigma} \rightarrow {\mathfrak M}(\sigma)\times X.
\end{equation}
By ~\cite[theorem II.1.7]{K}, to check that a point
$\zeta =\lt((C_v), (q_{f'}), (h_v)\rt)$ of $\Kbm{X}{\sigma}$ is in the
smooth locus of $(\pi,\text{ev}_f)$, it suffices to check that
$H^1(C,h^*T_X(-q_f)) = 0$.  For a general point in the image of
$\Kbmo{X}{\alpha}$, this follows from item $(3)$ of $\cB(X,\tau,f)$
along with lemma~\ref{lem-Breform2}.
\end{proof}

\begin{cor}~\label{cor-Bdefm}
Suppose $X\subset \PP^N$ is a smooth subvariety which satisfies
$\cB(X,\tau_1(e),f_1)$ for all $e=1,\dots,E$.  Let $\tau$ be an
$A$-graph with $E(\tau) \leq E$ and suppose that
$\alpha:\tau\rightarrow \sigma$ is a contraction.  The morphism
$\Kbmo{X}{\alpha}$ maps a general point of $\Kbmo{X}{\tau}$ to a
smooth point of $\Kbm{X}{\sigma}$.
\end{cor}

\begin{proof}
This is trivial if $\tau$ is empty.  Suppose $\tau$ is not empty and
let $v$ be a vertex of $\tau$.  Let $\tau\hookleftarrow \tau'$ be the
combinatorial morphism which attaches a new tail, $f$, at the vertex
$v$.  Let $\sigma\hookleftarrow \sigma'$ be the combinatorial morphism
which attaches a new tail, $f$, at the vertex of $\sigma$ which is the
image of $v$.  Let $\alpha':\tau'\rightarrow \sigma'$ be the
contraction which restricts to $\alpha$ and which maps $f$ to $f$.  
By proposition~\ref{prop-basic}, $\cB(X,\tau',f)$ holds.
By
proposition~\ref{prop-Bdefm}, $\Kbmo{X}{\alpha'}$ maps a general point
of $\Kbmo{X}{\tau'}$ to a point in the smooth locus of
$\text{ev}_f:\Kbm{X}{\sigma'}\rightarrow X$.  A point in the smooth
locus of $\text{ev}_f$ is a smooth point of $\Kbm{X}{\sigma'}$.  The
image of this point in $\Kbm{X}{\sigma}$ is also a smooth point.
Since $\Kbmo{X}{\tau'}$ surjects onto $\Kbmo{X}{\tau}$, a general
point of $\Kbmo{X}{\tau'}$ maps to a general point of
$\Kbmo{X}{\tau}$.  Thus $\Kbmo{X}{\alpha}$ maps a general point of
$\Kbmo{X}{\tau}$ to a smooth point of $\Kbm{X}{\sigma}$.  
\end{proof}

\textbf{Remark:}
Now suppose $\cB(X,\tau_1(e),f_1)$ holds for all $e=1,\dots,E$,
suppose that $\tau$ is a stable $A$-graph with $E(\tau)\leq E$ and
suppose that $\alpha:\tau\rightarrow \sigma$ is a contraction.  By
corollary~\ref{cor-Bdefm} and $(2)$ of proposition~\ref{prop-basic},
we conclude that there is a unique irreducible component $M(\alpha)$
of $\Kbm{X}{\sigma}$ which contains the image of $\Kbmo{X}{\alpha}$,
and $M(\alpha)$ is smooth of the expected dimension at a general
point.  So steps $(2)$ and $(3)$ of the strategy in the last section
are successful.  

\

Finally we give a simpler criterion for when $\cB(X,\tau,f)$ holds for all
$\tau$ with $E(\tau)\leq E$, where $E$ is some fixed integer, and also
reduce the number of components $M(\alpha)$ we have to deal with in
step $(4)$ of our strategy.

\begin{prop}~\label{prop-simplify}
Suppose that $X\subset \PP^N$ is a smooth subvariety satisfying
\begin{enumerate}
\item $\cB(X,\tau_1(1),f_1)$ holds,
\item $\fe(X,\tau_1(e),f_1)$ holds for $e=1,\dots,E$, and
\item $\Kbmo{X}{\tau_0(e)}$ is irreducible for $e=1,\dots,E$.
\end{enumerate}
Then for each stable $A$-graph $\tau$ with $E(\tau)\leq E$ and each
flag $f\in\text{Flag}(\tau)$, $\cB(X,\tau,f)$ holds and there is a
nice contraction $\alpha:\sigma\rightarrow \tau$ such that
$E(\sigma) \leq 1$ and such that $\Kbmo{X}{\alpha}$ maps the general
point of $\Kbmo{X}{\sigma}$ to a smooth point of $\Kbm{X}{\tau}$.
\end{prop}

\begin{proof}
It is easy to see that there is always a nice contraction
$\alpha:\sigma\rightarrow \tau$ such that $E(\sigma) \leq 1$.  By
proposition~\ref{prop-Bdefm}, we know that $\Kbmo{X}{\alpha}$ maps a
general point of $\Kbmo{X}{\sigma}$ to a smooth point of
$\Kbm{X}{\tau}$.  

\

By proposition~\ref{prop-basic}, to prove that $\cB(X,\tau,f)$ 
holds for all $\tau$ with $E(\tau)\leq E$, 
it suffices to prove that $\cB(X,\tau_1(e),f_1)$ holds for all
$e=1,\dots,E$.  
Let $\Kbm{X}{\tau_1(e)}'$ denote the normalization
of $\Kbm{X}{\tau_1(e)}$.  Consider the Stein factorization
$\Kbm{X}{\tau_1(e)}'\rightarrow Z\rightarrow X$ of $\text{ev}_{f_1}$.  By
proposition~\ref{prop-Bdefm}, there is a open, dense subset $U\subset
\Kbmo{X}{\sigma}$ such that $\Kbmo{X}{\alpha}$ maps $U$ into the
smooth locus of $\text{ev}_{f_1}$.  So $\Kbmo{X}{\alpha}|_U:U\rightarrow
\Kbm{X}{\tau_1(e)}$ factors through $\Kbm{X}{\tau_1(e)}'$.  Now consider the
image $V$ of $U$ in $Z$.  By proposition~\ref{prop-basic}, the general
fiber of $\text{ev}_{f_1}|_U:U\rightarrow X$ is geometrically
irreducible.  Therefore $V\rightarrow X$ is generically injective.

\

Let $\overline{V}\subset Z$ be the Zariski closure of $V$ with the
induced, reduced scheme structure.  Then $\overline{V}$ is an
irreducible stack and $\overline{V}\rightarrow X$ is surjective and
generically injective.  In particular, $\overline{V}$ is a nonempty
irreducible component of $Z$.  So the preimage of $V$ in
$\Kbm{X}{\tau_1(e)}$ is a nonempty irreducible component of
$\Kbm{X}{\tau_1(e)}$.  By corollary~\ref{cor-bbflat2}, every stratum in the
Behrend-Manin decomposition of $\Kbm{X}{\tau_1(e)}$ has the expected
dimension.  Thus $\Kbmo{X}{\tau_1(e)}$ is Zariski dense in
$\Kbm{X}{\tau_1(e)}$.  By assumption, $\Kbmo{X}{\tau_0(e)}$ is
irreducible.  Since $\Kbmo{X}{\tau_1(e)}\rightarrow
\Kbmo{X}{\tau_1(e)}$ is smooth with geometrically connected fibers,
also $\Kbmo{X}{\tau_1(e)}$ is irreducible.  Therefore
$\Kbm{X}{\tau_1(e)}$ is irreducible.  Therefore $\overline{V}=Z$ and
we conclude that the general fiber of $\Kbm{X}{\tau_1(e)}'\rightarrow
X$ is normal and geometrically connected, thus geometrically
irreducible.  It follows that the general fiber of
$\Kbmo{X}{\tau_1(e)}\rightarrow X$ is geometrically connected, i.e. we
have established item $(2)$ of the definition of $\cB(X,\tau_1(e),f_1)$.

\

To establish item $(3)$ of the definition of $\cB(X,\tau_1(e),f_1)$,
observe that the locus of points in $\Kbm{X}{\tau_1(e)}$ parametrizing
stable maps for which $h^*T_X$ is generated by global sections is an
open locus.  By lemma~\ref{lem-Breform} and the assumption
$\cB(X,\tau_1(1),f_1)$, for a general point of $\Kbmo{X}{\sigma}$, we
have that $h^*T_X$ is generated by global sections.  So this open set
intersects the general point of the image of $\Kbmo{X}{\alpha}$, so it
is nonempty.  Therefore it intersects $\Kbmo{X}{\tau}$ and item $(3)$
follows.  
\end{proof}

We summarize the results of this section for the case of complete
intersections in the following corollary.

\begin{cor}~\label{cor-Bsum}
Suppose that $X\subset \PP^N$ is a smooth complete intersection of
threshold degree $E(X)$ which satisfies:
\begin{enumerate}
\item $\cB(X,\tau_1(1),f_1)$ holds,
\item $\fe(X,\tau_1(e),f_1)$ holds for $e=1,\dots,E(X)$, and
\item $\Kbmo{X}{\tau_0(e)}$ is irreducible for $e=1,\dots,E(X)$.
\end{enumerate}
Then we have
\begin{enumerate}
\item For each basic $A$-graph $\tau$ and each flag
  $f\in\text{Flag}(\tau)$, $\cB(X,\tau,f)$ holds.
\item For each stable $A$-graph $\tau$ and each contraction
  $\alpha:\sigma\rightarrow \tau$ of a basic $A$-graph $\sigma$ to
  $\tau$, there is a unique irreducible component $M(\alpha)$ of
  $\Kbm{X}{\tau}$ which contains the image of $\Kbmo{X}{\alpha}$.
  Moreover $M(\alpha)$ is smooth of the expected
  dimension at a general point of the image of $\Kbmo{X}{\alpha}$.
\item $\Kbm{X}{\tau}$ is the union of the irreducible components
  $M(\alpha)$ as $\alpha:\sigma\rightarrow \tau$ ranges over nice
  contractions such that $E(\sigma)\leq 1$.  
\end{enumerate}
\end{cor}

\begin{proof}
The only new statement is $(3)$.  By theorem~\ref{thm-spec}, we know
that each irreducible component of $\Kbm{X}{\tau}$ is one of the
irreducible components $M(\alpha')$ for a nice contraction
$\alpha':\sigma'\rightarrow \tau$ with $\sigma$ a basic $A$-graph.  Of
course we can find a nice contraction $\beta:\sigma\rightarrow
\sigma'$ such that $E(\sigma) \leq 1$.  Let $\alpha:\sigma\rightarrow
\tau$ be the composition of $\beta$ and $\alpha'$.  Then $M(\alpha')$
is an irreducible component which contains the image of
$\Kbmo{X}{\alpha}$.  So $M(\alpha') = M(\alpha)$, i.e. we have proved
$(3)$.
\end{proof}

\section{Equating Irreducible Components}

Suppose that $X\subset \PP^N$ is a complete intersection which
satisfies the hypotheses of corollary~\ref{cor-Bsum}.  Then for each
stable $A$-graph $\tau$, we know that $\Kbm{X}{\tau}$ has the expected
dimension and is a union of irreducible components $M(\alpha)$ as
$\alpha:\sigma\rightarrow \tau$ ranges over nice contractions
with $E(\sigma)\leq 1$.  To prove that $\Kbm{X}{\tau}$ (and hence
$\Kbmo{X}{\tau}$) is irreducible, we are reduced to proving that the
irreducible components $M(\alpha)$ are all equal.

Suppose that $\cB(X,\tau_1(e),f_1)$ holds for all $e=1,\dots,E$ where
$E$ is some integer with $E\geq E(X)$. 
Fix a stable $A$-graph $\tau$ and
let $S_E(\tau)$ be the set of (isomorphism classes of) nice contractions
$\alpha:\sigma\rightarrow \tau$ with $E(\sigma)\leq E$.
Define a relation $\alpha \leq \alpha'$ if there exists a contraction
$\epsilon:\sigma \rightarrow \sigma'$ such that $\alpha=\alpha'\circ
\epsilon$.  If $\alpha\leq \alpha'$, then observe
$M(\alpha)=M(\alpha')$.  Form the equivalence relation $\cong$ on $S_E(\tau)$
generated by $\leq$.  Notice conclusion $(3)$ of
corollary~\ref{cor-Bsum} implies that every equivalence class contains
a contraction $\alpha:\sigma\rightarrow \tau$ such that $E(\sigma)\leq
1$.  Since $M(\alpha)=M(\alpha')$ if $\alpha\cong \alpha'$, we see
that the number of irreducible components of $\Kbm{X}{\tau}$ is
bounded by the number of equivalence classes of $\cong$ on $S_E(\tau)$.  So to
prove that $\Kbm{X}{\tau}$ is irreducible, it suffices to prove that
every two elements of $S_E(\tau)$ are equivalent.

\begin{defn}~\label{defn-thresh'}
Given $X\subset \PP^N$ a smooth complete intersection, define the
\emph{modified threshold degree} of $X$ to be $E'(X) = \max(E(X),2)$.
\end{defn}

\begin{prop}~\label{prop-abasic}
Suppose that $X\subset \PP^N$ is a smooth complete intersection such
that for $e=1,\dots,E'(X)$, we have $\cB(X,\tau_1(e),f_1)$ holds.
Then for each positive integer $e$, every two elements of
$S_{E'(X)}(\tau_0(e))$ are equivalent.  In particular
$\Kbm{X}{\tau_0(e)}$ is irreducible.
\end{prop}

\begin{proof}
Recall a connected tree $\tau$ is called an \emph{path} if $\tau$ has
precisely 
one or two leaves (so no vertex has valence greater than $2$).  
The number of vertices in a \emph{path} is the
\emph{diameter} of the path.  Given any connected tree $\tau$, the
\emph{diameter} of $\tau$, $\text{diam}(\tau)$, 
is defined to be the maximum diameter of a
subgraph which is a path.  If $\alpha:\sigma\rightarrow \tau_0(e)$ is a
nice contraction, then there are at most $e$ vertices in $\sigma$.  So
the diameter of $\sigma$ is at most $e$.  Moreover, there is a unique
contraction $\alpha_e:\sigma_e\rightarrow \tau_0(e)$ with
$\text{diam}(\sigma_e)=e$.  Here $\sigma_e$ is the $A$-graph whose
underlying graph is the path of length $e$, and for each vertex $v\in
\sigma_e$ we have $\beta(v)=1$.  

\

To prove that any two elements in $S_{E'(X)}(\tau_0(e))$ are
equivalent, it suffices to prove that any two nice contractions
$\alpha:\sigma\rightarrow \tau$ with $E(\sigma)=1$ are equivalent.  We
will prove that for each such $\alpha:\sigma\rightarrow \tau$ with
$\text{diam}(\sigma) < e$, there is a nice contraction
$\alpha':\sigma'\rightarrow \tau$ such that $\sigma\cong \sigma'$,
$E(\sigma')=1$ and $\text{diam}(\sigma') \geq \text{diam}(\sigma)$.
From this it follows by induction that all such contractions are
equivalent to $\alpha_e:\sigma_e\rightarrow \tau$.  

\

Suppose that $\alpha:\sigma\rightarrow \tau$ is a nice contraction
with $E(\sigma)=1$ and $\text{diam}(\sigma) < e$.  Let
$\gamma\hookrightarrow \sigma$ be a subgraph which is a path such that
$\text{diam}(\gamma)=\text{diam}(\sigma)$.  Since $\gamma$ does not
equal $\sigma$, there exists a vertex $v_1$ of $\gamma$ such that the
valence of $v_1$ is at least $3$.  Let ${f_1,f_2}$ be an edge of
$\sigma$ not contained in $\gamma$ such that $\partial f_1 = v_1$.
Let $v_2=\partial f_2$.  Form the nice contraction $\epsilon:\sigma\rightarrow
\rho$ which contracts $v_1$ and $v_2$ to a common vertex $v$ of
$\sigma$ with $\beta(v)=2$.  The nice contraction
$\alpha:\sigma\rightarrow \tau_0(e)$ factors through a nice contraction
$\alpha_\rho:\rho\rightarrow \tau_0(e)$.

\

The image of $\gamma$ in $\rho$ is a path
$\gamma_\rho$ which contains $v$.  Now let $\gamma'\rightarrow \gamma_\rho$
be a contraction of a path of length $\text{diam}(\gamma)+1$ which
contracts two adjacent vertices $w_1$ and $w_2$ to $v$ (where
$\beta(w_1)=\beta(w_2)=1$).  There is a unique nice contraction
$\epsilon':\sigma'\rightarrow \rho$ such that $\gamma'$ is a path in
$\sigma'$, such that the restriction of $\epsilon'$ to $\gamma'$ is just
$\gamma'\rightarrow \rho$, such that every flag of $v$ not contained
in $\gamma_\rho$ is the image of a flag of $w_1$, and which is an
isomorphism from $\sigma'-\gamma'\rightarrow \rho-\gamma_\rho$.
Define $\alpha'=\alpha_\rho\circ \epsilon$.  Then
$\alpha':\sigma'\rightarrow \tau_0(e)$ is a nice contraction,
$E(\sigma')=1$, $\alpha\cong \alpha'$ and
$\text{diam}(\sigma')=\text{diam}(\sigma)+1$.  This proves the claim.

\

So by induction on the diameter of $\sigma$, every element of
$S_{E'(X)}(\tau_0(e))$ is equivalent to $\alpha_e:\sigma_e\rightarrow
\tau_0(e)$.  In particular $\Kbm{X}{\tau_0(e)}=\Kgnb{0}{0}{X}{e}$ is
irreducible.  
\end{proof}

\begin{cor}~\label{cor-abasic}
With the same hypotheses as in proposition~\ref{prop-abasic}, for each
stable $A$-graph $\tau$ we have
\begin{enumerate}
\item $\Kbm{X}{\tau}$ is an integral, local complete intersection
  stack of the expected dimension $\text{dim}(X,\tau)$, and
  $\Kbmo{X}{\tau}$ is the unique dense stratum in the Behrend-Manin
  decomposition.  
\item For each flag $f\in \text{Flag}(\tau)$, $\cB(X,\tau,f)$ holds.
\item For each contraction $\alpha:\sigma\rightarrow \tau$,
  $\Kbm{X}{\tau}$ is smooth at the general point of the image of
  $\Kbmo{X}{\alpha}: \Kbmo{X}{\sigma} \rightarrow \Kbm{X}{\tau}$.
\end{enumerate}
\end{cor}

\begin{proof}
By corollary~\ref{cor-bbflat}, $\fe(X,\tau_1(e),f_1)$ holds for all
integers $e>0$.  By assumption, $\cB(X,\tau_1(1),f_1)$ holds.  And by
proposition~\ref{prop-abasic}, $\Kbmo{X}{\tau_0(e)}$ is irreducible
for each integer $e>0$.  Thus by proposition~\ref{prop-simplify}, for
every stable $A$-graph $\tau$ and every flag $f\in\text{Flag}(\tau)$,
we have that $\cB(X,\tau,f)$ holds.  This establishes $(2)$.  

\

As in the proof of
proposition~\ref{prop-basic}, $(2)$ implies that for every stable
$A$-graph $\tau$, $\Kbmo{X}{\tau}$ is irreducible of the expected
dimension.  
By a parameter count, we conclude that $\Kbmo{X}{\tau}$ is
the unique dense stratum of the Behrend-Manin decomposition of
$\Kbm{X}{\tau}$.  So $\Kbm{X}{\tau}$ is also irreducible of the
expected dimension, and generically smooth.  So $\cD(X,\tau)$ holds.
By lemma~\ref{lem-propeq}, we conclude that $\lci(X,\tau)$ holds,
i.e. $\Kbm{X}{\tau}$ is a local complete intersection stack.  Since it
is generically smooth, and thus generically reduced, it is reduced.
So $\Kbm{X}{\tau}$ is an integral, local complete intersection stack
of the expected dimension $\text{dim}(X,\tau)$ and $\Kbmo{X}{\tau}$ is
the unique dense stratum in the Behrend-Manin decomposition.  This
establishes $(1)$.

\

Finally $(3)$ follows from $(1)$ and corollary~\ref{cor-Bdefm}.
\end{proof}

Finally, we prove that a general hypersurface $X\subset \PP^n$ of
degree $d<\frac{n+1}{2}$ satisfies the hypotheses of
proposition~\ref{prop-abasic}.  

\begin{prop}~\label{prop-final}
Suppose $n>2$, $d < \frac{n+1}{2}$ and suppose $X\subset \PP^n$ is a
hypersurface of degree $d$, so $E'(X)=2$.  If $\cB(X,\tau_1(1),f_1)$
holds for $X$ (recall from lemma~\ref{lem-Blines} that
$\cB(X,\tau_1(1),f_1)$ holds for a general $X$) then also
$\cB(X,\tau_1(2),f_1)$ holds.  For such an $X$, the results of
corollary~\ref{cor-abasic} hold.
\end{prop}

\begin{proof}
By
corollary~\ref{cor-bbflat}, $\lci(X,\tau_1(2),f_1)$ holds.  
Since $\Kbmo{X}{\tau_0(2)}$ is the unique dense stratum in the
Behrend-Manin decomposition, to prove that $\Kbmo{X}{\tau_0(2)}$ is
irreducible, it is equivalent to prove that $\Kbm{X}{\tau_0(2)}$ is
irreducible.  
To see that $\Kbm{X}{\tau_0(2)}$ is irreducible, observe by
theorem~\ref{thm-spec} that every
irreducible component of $\Kbm{X}{\tau}$ is of the form $M(\alpha)$
for a nice contraction $\alpha:\sigma\rightarrow \tau_0(2)$ with
$E(\alpha)=1$. But there is a unique such contraction, namely
$\alpha_2:\sigma_2\rightarrow \tau_0(2)$.  So $\Kbmo{X}{\tau_0(2)}$ is
irreducible.  

\

Now by property~\ref{prop-simplify}, $\cB(X,\tau_1(2),f_1)$ holds.  
\end{proof}

\bibliography{my}

\begin{thebibliography}{10}

\bibitem{B}
K.~Behrend.
\newblock Gromov-{W}itten invariants in algebraic geometry.
\newblock {\em Inventiones Mathematica}, 127:601--617, 1997.
\newblock ar{X}iv:math.AG/9601011.

\bibitem{BM}
K.~Behrend and Y.~Manin.
\newblock Stacks of stable maps and {G}romov-{W}itten invariants.
\newblock {\em Duke Math Journal}, 85:1--60, 1996.
\newblock ar{X}iv:math.AG/9506023.

\bibitem{FP}
W.~Fulton and R.~Pandaharipande.
\newblock Notes on stable maps and quantum cohomology.
\newblock In {\em Algebraic geometry -- Santa Cruz 1995}, pages 45--96.
  American Mathematical Society, 1995.
\newblock ar{X}iv:math.AG/9608011.

\bibitem{SGA2}
A.~Grothendieck.
\newblock {\em Cohomologie locale des faisceaux coh{\'e}rents ... (SGA 2)},
  volume~2 of {\em Advanced Studies in Pure Mathematics}.
\newblock North-{H}olland {P}ublishing {C}o., 1968.

\bibitem{HS2}
J.~Harris and J.~Starr.
\newblock Rational curves on hypersurfaces of low degree, {II}.
\newblock in preparation.

\bibitem{H}
R.~Hartshorne.
\newblock {\em Algebraic Geometry}, volume~52 of {\em Graduate Texts in
  Mathematics}.
\newblock Springer-Verlag, 1977.

\bibitem{KP}
B.~Kim and R.~Pandharipande.
\newblock The connectedness of the moduli space of maps to homogeneous spaces.
\newblock In {\em Symplectic geometry and mirror symmetry}, chapter~5. World
  Scientific, 2001.
\newblock ar{X}iv:math.AG/0003168.

\bibitem{K}
J.~Koll\'ar.
\newblock {\em Rational Curves on Algebraic Varieties}, volume~32 of {\em
  Ergebnisse der Mathematik und ihrer Grenzgebiete, 3. Folge}.
\newblock Springer-Verlag, 1996.

\bibitem{Ma}
H.~Matsumura.
\newblock {\em Commutative ring theory}, volume~8 of {\em Cambridge studies in
  advanced mathematics}.
\newblock Cambridge University Press, 1986.

\bibitem{AVar}
D.~Mumford.
\newblock {\em Abelian varieties}, volume~5 of {\em Tata Institute of
  Fundamental Research Studies in Mathematics}.
\newblock Oxford University Press, 1970.

\end{thebibliography}
\bibliographystyle{abbrv}

\end{document}